\documentclass[12pt]{article}
\usepackage{amsmath,amsthm,amssymb,amscd}
\usepackage{latexsym}
\newtheorem{thm}{Theorem}[section]
 \newtheorem{cor}[thm]{Corollary}
 \newtheorem{lem}[thm]{Lemma}
 \newtheorem{prop}[thm]{Proposition}
 \theoremstyle{definition}
 
 \theoremstyle{remark}

  \numberwithin{equation}{section}

\makeatletter
\newcommand*\owedge{\mathpalette\@owedge\relax}
\newcommand*\@owedge[1]{%
  \mathbin{%
    \ooalign{%
      $#1\m@th\bigcirc$\cr
      \hidewidth$#1\m@th\wedge$\hidewidth\cr
    }%
  }%
}
\makeatother

\begin{document}

\title
{Classification of  compact manifolds with positive isotropic curvature}

\author{ Hong Huang}
\date{}
\maketitle

\begin{abstract}
 We show the following result: Let $(M,g_0)$ be a compact connected
manifold of dimension $n\geq 12$ with  positive isotropic curvature. Then $M$ is diffeomorphic  to a spherical space form, or a  quotient manifold of $\mathbb{S}^{n-1}\times \mathbb{R}$ by a cocompact discrete subgroup of the isometry group of the round cylinder $\mathbb{S}^{n-1}\times \mathbb{R}$, or a connected sum of a finite number of such manifolds.
This extends   previous works of Brendle and Chen-Tang-Zhu,  and improves a work of Huang.  The proof uses Ricci flow with surgery on compact orbifolds, with the help of  the ambient isotopy uniqueness of closed tubular neighborhoods of  an isolated singular point in an orbifold.


\end{abstract}


\section {Introduction}

In our previous work \cite{Hu22}, based on Brendle's remarkable paper \cite{B19}, we get a rough result on the classification  of compact manifolds with positive isotropic curvature of dimension $n\geq 12$, from which a conjecture of Schoen and a conjecture of Gromov (and of Fraser), both in dimensions $n\geq 12$, follow. Now we can improve the main result in \cite{Hu22}, and get a more precise topological classification of  compact manifolds with positive isotropic curvature of dimension $n\geq 12$.

 \begin{thm} \label{thm 1.1} Let $(M,g_0)$ be a compact connected
manifold of dimension $n\geq 12$ with  positive isotropic curvature. Then $M$ is diffeomorphic  to a spherical space form, or a  quotient manifold of $\mathbb{S}^{n-1}\times \mathbb{R}$ by a cocompact discrete subgroup of the isometry group of the round cylinder $\mathbb{S}^{n-1}\times \mathbb{R}$, or a connected sum of a finite number of such manifolds.
\end{thm}

 \noindent (For the classification of spherical space forms see \cite{Wo}. Of course, the cocompact discrete subgroup of the isometry group of the round cylinder $\mathbb{S}^{n-1}\times \mathbb{R}$ referred to in the above theorem should act freely on $\mathbb{S}^{n-1}\times \mathbb{R}$; moreover, it is well known that $\text{Isom}(\mathbb{S}^{n-1}\times \mathbb{R})\cong \text{Isom}(\mathbb{S}^{n-1}) \times \text{Isom}(\mathbb{R})$ for $n\geq 2$ (cf. for example, the proof of Corollary 6.2 in \cite{CG}). For the investigation of compact quotient manifolds of $\mathbb{S}^{n-1}\times \mathbb{R}$ by standard isometries one can consult for example, p. 457 of \cite{Sc}, p. 827 of \cite{CH}, Theorem 1.1 in \cite{GSM},  and \cite{MS}. By  \cite{MW} the converse of Theorem \ref{thm 1.1} is also true.) This extends the Main Theorem in Chen-Tang-Zhu \cite{CTZ}, which deals with the 4 dimensional case, and Theorem 1.4 in Brendle \cite{B19}, which
deals with the case that the $n$-manifolds contain no nontrivial incompressible $(n-1)$-dimensional space forms, where $n \geq 12$.  The result is analogous to the classification of compact 3-manifolds with positive scalar curvature by Perelman et al.

As an application we partially verify a conjecture in \cite{Hu22}.

\begin{cor} \label{cor 1.2} \ \  Let $n\geq 12$ be even, $\Gamma$ be a  finite subgroup of $O(n)$ acting freely on $\mathbb{S}^{n-1}$,  and   $M$ be the total space of an orientable $\mathbb{S}^{n-1}/\Gamma$-bundle over $\mathbb{S}^1$. Suppose that $M$ admits a metric with positive isotropic curvature. Then  $M$ admits $\mathbb{S}^{n-1} \times \mathbb{R}$-geometry.
\end{cor}

\noindent (By $M$ admits $\mathbb{S}^{n-1} \times \mathbb{R}$-geometry we mean that $M$ admits a Riemannian metric which is locally isometric to the round cylinder $\mathbb{S}^{n-1} \times \mathbb{R}$; in this case $M$ is diffeomorphic to a quotient of the round cylinder $\mathbb{S}^{n-1} \times \mathbb{R}$ by some isometries.)

The following result is a simple consequence of Theorem \ref{thm 1.1}, which strengthens and extends the Main Theorem in Micallef-Moore \cite{MM} in the case of $n \geq 12$.

\begin{cor} \label{cor 1.3} \ \  Let $(M,g_0)$ be a compact connected
	manifold of dimension $n\geq 12$ with  positive isotropic curvature. Suppose that the fundamental group of $M$ is finite. Then $M$ is diffeomorphic  to a spherical space form.
\end{cor}

\noindent In particular, since  any compact manifold which admits a non-collapsed, ancient Ricci flow must have finite fundamental group, cf.  \cite{Ba}, it follows that any compact (connected) manifold of dimension $n\geq 12$  which admits a non-collapsed, ancient Ricci flow with positive isotropic curvature is diffeomorphic to a spherical space form. (Recall that by combining Perelman's result that  ``there is no shrinking breathers other than gradient
solitons'' on a closed manifold  (see p.9 of \cite{P1}, cf. also Theorem 1.3 in \cite{N}) with  Theorem 1.4 in  \cite{Na} one sees that any compact shrinking Ricci soliton of dimension $n\geq 12$  with positive isotropic curvature is isometric  to a spherical space form.)

Recall that Micallef and Moore \cite{MM}  showed that for a  compact
manifold $M$ of dimension $n\geq 4$ with  positive isotropic curvature, we have the homotopy groups $\pi_i(M)=0$ for $2\leq i \leq \lfloor \frac{n}{2} \rfloor$; compare Theorem 4.6(ii) in  Schoen and Yau \cite{SY}.  The following corollary extends  Micallef-Moore's result when $n\geq 12$.

\begin{cor} \label{cor 1.4} \ \  Let $M$ be a compact smooth
manifold  of dimension $n\geq 12$ which admits a metric with  positive isotropic curvature. Then  $\pi_i(M)=0$ for $2\leq i \leq n-2$.
\end{cor}

Note that very recently Chow and Wang \cite{CW} showed that Corollary 1.4 is also true when $n=5$ and $n=6$.

In \cite{MW} Micallef-Wang  showed that a compact manifold $M$ of even dimension $n=2m \geq 4$ with positive isotropic curvature has the second Betti number $b_2(M)=0$. The following corollary extends  their result when $n\geq 12$ (for our result the dimension $n$ needs not be even).

\begin{cor} \label{cor 1.5} \ \  Let $M$ be a compact smooth
manifold  of dimension $n\geq 12$ which admits a metric with  positive isotropic curvature. Then  $b_i(M)=0$ for $2\leq i \leq n-2$.
\end{cor}

As in the proof of Corollary 2 in \cite{CTZ}, by using Corollary 5 in \cite{SY79} we also have the following

\begin{cor} \label{cor 1.6} \ \: Let $M$ be a compact smooth
	manifold of dimension $n\geq 12$. If $M$  admits a  metric with  positive isotropic curvature, then $M$ also admits a locally conformally flat  metric with  positive scalar curvature.
\end{cor}

\noindent (But the converse is not true; see \cite{MW}.)

In the process of proof of Theorem \ref{thm 1.1} we get a slightly more general result.

\begin{thm} \label{thm 1.7} \ \  Let $(\mathcal{O},g_0)$ be a compact connected
orbifold of dimension $n\geq 12$ with at most isolated singularities and with  positive isotropic curvature. Then $\mathcal{O}$ is diffeomorphic to  an orbifold connected sum of a finite number of spherical orbifolds with at most isolated singularities.
\end{thm}

\noindent (By an obvious extension of Theorem 1.1 in \cite{MW} to the case of orbifold connected sum (cf. also \cite{CH}), a certain version of the converse of Theorem \ref{thm 1.7} is also true.) This improves Theorem 1.3 in \cite{Hu22}. Compare Theorem 2.1 in \cite{CTZ}.  We emphasize that here our notion of orbifold connected sums (see \cite{Hu22}) is stronger than that defined on p. 48 of \cite{CTZ}.  (As in \cite{CTZ} we also allow self-connected sums; but of course our notion of self-connected sums is also stronger than that given in \cite{CTZ}.)

As in \cite{CTZ}, \cite{B19}, and \cite{Hu22}, the main tool in our proof is  Ricci flow with surgery, aided by some topological arguments.
In the 4-dimensional case the topological part of the  proof of the Main Theorem in Chen-Tang-Zhu \cite{CTZ}  used a result of McCullough \cite{Mc}, which says that any diffeomorphism of a 3-dimensional spherical space form is isotopic to an isometry.  Similar statement for $(n-1)$-dimensional spherical space forms with $n \geq 5$ is not true in general.

One of our key observations is that for the topological part in the higher dimensional case, it suffices to use  the ambient isotopy uniqueness of closed tubular neighborhoods of  an isolated singular point in an orbifold.  More precisely, using this uniqueness and a (rough) classification of spherical orbifolds with (nonempty) isolated singularities (due to  \cite{CTZ} in the even ($\geq 4$) dimensional case and to \cite{Hu22} in the odd ($\geq 3$) dimensional case) one can prove  fine properties of orbifold ancient $\kappa$-solutions and standard solutions when they contain caps with ends of type $\mathbb{S}^{n-1}/\Gamma \times (-1,1)$ with $\Gamma$ nontrivial, extending  what were done in \cite{B19} and  \cite{Hu22}  for caps with ends of type $\mathbb{S}^{n-1} \times (-1,1)$,   and extend/strengthen the  definition of $\varepsilon$-caps in  \cite{B19} and   \cite{Hu22}. Then we can recognize the topology of components in the surgery process which are covered by canonical neighborhoods. Moreover, using the ambient isotopy uniqueness of closed tubular neighborhoods of  an isolated singular point again, one can determine the topology of a connected sum of at most two spherical orbifolds with at most isolated singularities.   (It is clear from our arguments in this paper that the use of McCullough's result in the proof of the Main Theorem in Chen-Tang-Zhu \cite{CTZ} can be avoided; moreover, to get a similar classification of compact manifolds of dimension $5\leq n \leq 11$ with positive isotropic curvature, the only thing that remains to do is to show  a suitable curvature pinching estimate, as Hamilton did in dimension 4 in \cite{H97} and Brendle did in dimensions $n\geq 12$ in \cite{B19}.)

We emphasize that since in a local chart (model) a local lift of an orbifold diffeomorphism is only weakly equivariant (that is, equivariant with respect to a group isomorphism not necessarily the identity map; cf. p. 31 of \cite{HKMR}), in our situation in general it does not suffice to simply apply the equivariant tubular neighborhood theorem  in \cite{B72}, which requires (strongly) equivariant map, in local models. But note that in this paper we only need to use the orbifold tubular neighborhood theorem (Proposition \ref{prop 2.1}) in the case that the ambient orbifold is a spherical orbifold with at most isolated singularities. In this case the orbifold tubular neighborhood theorem is implied by  a variant of equivariant tubular neighborhood theorem for an ambient manifold diffeomorphic to either $\mathbb{S}^n $ or $\mathbb{R}^n $ and with a finite group action, where  one should allow weakly equivariant embeddings in the definition of invariant tubular neighborhoods (compare \cite{B72}). It turns out that one can prove such a variant of equivariant tubular neighborhood theorem by adapting some ideas from the existing literature.  In fact we'll show a slightly more general orbifold version of the tubular neighborhood theorem, whose precise statement and detailed proof are not easy to find in the existing literature. Here we also realize that it is very convenient to use possibly noneffective orbifolds and orbifold vector bundles over possibly noneffective orbifolds.

 Recently,  Zhengnan Chen \cite{Ch} /respectively Jae Ho Cho \cite{Ch} extended Brendle's /respectively Brendle's and Chen's curvature pinching estimates to the case of $n\geq 9$ /respectively $n\geq 8$. So combined with their estimates, our arguments imply that all our results above also hold true for $8 \leq n \leq 11$. Chow and Wang \cite{CW} showed Schoen's conjecture in the cases of $n=5$ and $n=6$.

 In Section 2 we show the ambient isotopy uniqueness of closed tubular neighborhoods of  certain special suborbifolds in an orbifold, define $\varepsilon$-caps,
and describe the canonical neighborhood structure of orbifold ancient
 $\kappa$-solutions and orbifold standard solutions. In Section 3 we  construct  an $(r, \delta)$-surgical solution to the Ricci flow starting with a compact, connected
Riemannian orbifold of dimension $n\geq 12$ with at most isolated singularities and with  positive isotropic curvature. As applications we prove
Theorems \ref{thm 1.7} and \ref{thm 1.1},    Corollaries \ref{cor 1.2}, \ref{cor 1.3}, \ref{cor 1.4} and \ref{cor 1.5}, and  derive a topological classification of compact Riemannian orbifold of dimension $n\geq 12$ with at most isolated singularities and with positive isotropic curvature. We will follow the notation and conventions in \cite{Hu22}  in general. In particular, as in \cite{Hu22}, when necessary we'll assume that the smooth maps between orbifolds that appear are smooth complete orbifold maps in the sense of \cite{BB12} and \cite{BB15}  to avoid any possible ambiguity.

\section{$\varepsilon$-caps and ancient $\kappa$-solutions}

In general we  use   $\Gamma$  to denote a  finite subgroup of $O(n)$ acting freely on $\mathbb{S}^{n-1}$.

\vspace *{0.2cm}

 To describe the structure of orbifold ancient $\kappa$-solutions and orbifold standard solutions we need  notions on necks and caps. Let $\mathcal{O}$ be an $n$-dimensional orbifold,  $U$  an open subset in $\mathcal{O}$ (sometimes when we write $\mathcal{O}$ we really mean its underlying space $|\mathcal{O}|$).
If $U$ is  diffeomorphic to $\mathbb{S}^{n-1}//\Gamma \times \mathbb{R}$, we call it  a topological neck.

 Suppose further that the spherical space form $\mathbb{S}^{n-1}/\Gamma$  admits  an isometric involution $\sigma$ with at most isolated fixed points, let $\hat{\sigma}$ be the
involution on the manifold $\mathbb{S}^{n-1}/\Gamma \times \mathbb{R}$
defined by $\hat{\sigma}(x,s)=(\sigma(x),-s)$ for $x\in
\mathbb{S}^{n-1}/\Gamma $ and $s\in \mathbb{R}$, consider the quotient orbifold $(\mathbb{S}^{n-1}/\Gamma \times
\mathbb{R})//\langle \hat{\sigma}\rangle$, which has at most  isolated singularities. We also denote this (Riemannian)
orbifold by $\mathbb{S}^{n-1}/\Gamma \times_{\mathbb{Z}_2}  \mathbb{R}$; but beware that this notation has some ambiguity.  Similarly we can define $\mathbb{S}^{n-1}/\Gamma \times_{\mathbb{Z}_2}  [-1,1]$. Note that we can consider $\Gamma$ and $\hat{\sigma}$ as isometries of $\mathbb{S}^n$ in a natural way, by lifting $\sigma$ to an isometry of  $\mathbb{S}^{n-1}$ (which is always possible) and viewing $\mathbb{S}^n$ as a suspension of
$\mathbb{S}^{n-1}$. We'll use the same notation for these isometries of $\mathbb{S}^n$.
The $(\mathbb{S}^{n-1}/\Gamma \times
\mathbb{R})//\langle \hat{\sigma}\rangle$  above is a smooth manifold if and only if $\sigma$ has no any fixed points in $\mathbb{S}^{n-1}/\Gamma$; if  this is the case, it is diffeomorphic to  $\mathbb{S}^n// \langle\Gamma, \hat{\sigma} \rangle \setminus \bar{B}$,  where $B$ is a small, open metric ball  centered at  the point coming from the north and  south poles of $\mathbb{S}^n$.  Note that when $n$ is odd, $(\mathbb{S}^{n-1}/\Gamma \times
\mathbb{R})//\langle \hat{\sigma}\rangle$  is a smooth manifold if and only if $\Gamma$ is  trivial   and  $\langle\sigma \rangle$   acts on $\mathbb{S}^{n-1}$  antipodally; see the proof of Proposition 2.5 in \cite{Hu22}.  Below when we write $\mathbb{S}^n// \langle\Gamma, \hat{\sigma} \rangle$, we  always assume that
the involution $\sigma$ has no any fixed points in $\mathbb{S}^{n-1}/\Gamma$.  (For the investigation of spherical space forms which admit free isometric involutions one can consult for example, \cite{GSM} and \cite{MS}.)

If $\mathbb{S}^{n-1}/\Gamma \times_{\mathbb{Z}_2}  \mathbb{R}$ has nonempty isolated singularities,  by Proposition 2.5 in \cite{Hu22}, it must be diffeomorphic to $\mathbb{S}^n// (x,\pm x') \setminus \bar{B}$, where, as in \cite{CTZ},   $\mathbb{S}^n// (x,\pm x')$ denotes the quotient orbifold   $\mathbb{S}^n// \langle\iota \rangle$ with $\iota$  the isometric involution on  $\mathbb{S}^n$ given by
 $(x_1,x_2,\cdot\cdot\cdot,x_{n+1}) \mapsto
(x_1,-x_2,\cdot\cdot\cdot,-x_{n+1})$, which  has two  singular points, and $B$ is a  small, open  metric ball around  a regular point in $\mathbb{S}^n// (x,\pm x')$ such that the closure  $\bar{B}$  of $B$  is disjoint from the two singular points and  diffeomorphic to $D^n$.  We'll choose the regular point to be the image of the north pole under the natural projection $\mathbb{S}^n \rightarrow  \mathbb{S}^n// (x,\pm x')$.

   If an open subset $U$ in an $n$-dimensional orbifold $\mathcal{O}$ is  diffeomorphic to  $\mathbb{R}^n//\Gamma$ (here $\Gamma$ is a finite subgroup of $O(n)$ acting freely on $\mathbb{S}^{n-1}$) or  $\mathbb{S}^{n-1}//\Gamma \times_{\mathbb{Z}_2}  \mathbb{R}$  (in this case the spherical space form $\mathbb{S}^{n-1}/\Gamma$ should  admit  an isometric involution $\sigma$ with at most isolated fixed points) defined above, we call it a topological cap.

\vspace*{0.2cm}

Let $\mathcal{P}$ be an  orbifold with or without boundary,  $\mathcal{O}$ be an orbifold without boundary. If $f: \mathcal{P} \rightarrow \mathcal{O}$ is an immersion (as defined on p. 107-108 of \cite{KL}) such that $|f|$ maps $|\mathcal{P}|$ homeomorphically onto its image in  $|\mathcal{O}|$, we call  $f$  an embedding and  $\mathcal{P}$ an embedded suborbifold of $\mathcal{O}$ (but beware that  the meaning of the word “embedded" here is not the same as  that on p. 108 of \cite{KL}). Let  $f, g: \mathcal{P}\rightarrow \mathcal{O}$ be two embeddings. Suppose that $h: \mathcal{P} \times \mathbb{R} \rightarrow \mathcal{O}$ is a smooth map,  such that for each $t\in \mathbb{R}$, $h_t:=h(\cdot, t)$ is an embedding, $h_0=f$, and $h_1=g$,  then we call $h$ an isotopy between the two embeddings $f$ and $g$ (or simply an isotopy of $f$); compare p. 206 in \cite{Mu} in the manifold case. Note that here we use  $\mathcal{P} \times \mathbb{R}$ instead of $\mathcal{P} \times [0,1]$, as when $\mathcal{P}$ has boundary, $\mathcal{P} \times [0,1]$ is an orbifold with corner, which we avoid considering; compare p. 64 in \cite{Mu}. If an isotopy $h$
between two embeddings $f$ and $g$ satisfies that for some small $\epsilon >0$, $h_t=f$ when $t \leq \epsilon$, and $h_t=g$ when  $t \geq 1-\epsilon$, we call $h$ a technical isotopy, following Definition 9.2 in \cite{BJ}.
A smooth map $H: \mathcal{O} \times \mathbb{R} \rightarrow \mathcal{O}$ with each $H_t$ being a diffeomorphism of $\mathcal{O}$  and $H_0=\text{Id}$ is called a diffeotopy of $\mathcal{O}$; compare p. 207 in \cite{Mu} in the manifold case. The support of a diffeotopy  $H$ of an orbifold $\mathcal{O}$ is the closure of the set of points $p \in \mathcal{O}$ such that $H_t(p) \neq p$ for some $t \in \mathbb{R}$.

Let $\mathcal{O}$ and $\mathcal{O}'$ be two orbifolds, and  $f: \mathcal{O}' \rightarrow \mathcal{O}$ be an embedding. We identify $\mathcal{O}'$ with its image in $\mathcal{O}$. If we assume further that  $\mathcal{O}'$ is a full suborbifold of $\mathcal{O}$ in the sense of Definition 5 in \cite{BB15} (see also \cite{We}; compare Definition 1.14 in \cite{D94}), wel call such an immersion $f$ (or embedded suborbifold $\mathcal{O}'$) an embedded full suborbifold (see also \cite{We} and Definition 2.5 in \cite{LZ}). There is also an obvious extension of the notion of embedded full suborbifold to the case of embedded full suborbifold with boundary, in which the orbifold $\mathcal{O}'$ above can have boundary, correspondingly the suborbifold $\mathcal{P}$  and the submanifold $\widetilde{V}_x$ in Definition 4 in \cite{BB15} can have boundary; consult also Sections 2.1.1 and 2.1.3 in \cite{BMP}. (Note that any point in $|\mathcal{O}|$ may be viewed as a 0-dimensional embedded full suborbifold of $\mathcal{O}$.)

Now we assume that $\mathcal{O}'$ is an embedded full suborbifold in $\mathcal{O}$ and view it as a possibly noneffective orbifold with local groups (whose actions on local models may  be noneffective) at all $p\in |\mathcal{O}'| \subset |\mathcal{O}|$ the same as the ones in the ambient orbifold $\mathcal{O}$. We still denote this possibly noneffective orbifold by $\mathcal{O}'$. Let $\xi$ be  an orbifold  vector bundle  over the possibly noneffective orbifold $\mathcal{O}'$ such that there is an embedding $\varphi: E(\xi) \rightarrow \mathcal{O}$ with $|\varphi|(|E(\xi)|)$ being an open neighborhood of $|\mathcal{O}'|$ in $|\mathcal{O}|$, and with the restriction of $\varphi$ to the zero-section of $\xi$ coinciding with the inclusion of $\mathcal{O}'$ in $\mathcal{O}$ (and the corresponding homomorphisms between local groups being isomophisms (but not necessarily the identity maps)), where $E(\xi)$ is  the total space of the bundle $\xi$, and is required  to be an effective orbifold. (The definition of an orbifold vector bundle over a possibly noneffective orbifold can be easily adapted from that of a vector orbi-bundle in Section 1.2.1 of \cite{GMS}; in particular, we extend Definition 1.45 of \cite{GMS} by allowing $\Gamma_2$ there to act noneffectively on $\hat{U}_2$.)  We call $\varphi$ an open tubular neighborhood of $\mathcal{O}'$ in $\mathcal{O}$.  For an embedded full suborbifold $\mathcal{O}'$ (without boundary) of $\mathcal{O}$ as above, such an orbifold  vector bundle  $\xi$ always exists; for example one can choose $\xi$  to be the normal bundle of the suborbifold $\mathcal{O}'$ in $\mathcal{O}$. (Compare Theorem 1.77 in \cite{GMS}.) Here we define the normal bundle of the embedded full suborbifold $\mathcal{O}'$ in $\mathcal{O}$ to be the orbifold vector bundle $i^*T\mathcal{O}/T\mathcal{O}'$ over the possibly noneffective orbifold $\mathcal{O}'$, where $i: \mathcal{O}' \rightarrow \mathcal{O}$ is the inclusion map  (from the possibly noneffective orbifold $\mathcal{O}'$ to the ambient orbifold $\mathcal{O}$) with all the corresponding  homomorphisms between local groups being the identity maps; consult p. 35-37 of \cite{ALR}, Remark 4.4.12c in \cite{CR}, and p. 304 of \cite{D94}. (The definition of the pullback of the tangent bundle used here can be obviously adapted from Definition 6.2 in \cite{BB13}; the only difference is that for any point $x$ in the underlying space of the domain orbifold, the local group $\Gamma_x$ (using the same notation as in \cite{BB13}) may act noneffectivily on the local model of the domain orbifold around $x$.) One can show this by adapting the proof of Theorem 12.11 in \cite{BJ} and Theorem 7.1.5 in \cite{Mu} to the orbifold case; consult also Theorem 2.2 in Chapter VI of \cite{B72}, where the closeness condition imposed on the invariant submanifold $A$ can be removed, Proposition 2.1.2 in \cite{D90} and Proposition 1.15 in \cite{D94}, whose proof uses Theorem 2.2 in Chapter VI of \cite{B72} and the orthonormal frame bundle of $\mathcal{O}$ (after introducing a Riemannian metric on $\mathcal{O}$), Theorem 7.3 in \cite{Ka}, and Theorem 2.10 in \cite{PW}.

\vspace*{0.2cm}

\noindent {\bf Example}.   Let $\mathcal{O}'$  be a point (a 0-dimensional connected manifold), $\mathcal{O}$ be an orbifold of dimension $n \geq 1$. Suppose $q \in |\mathcal{O}|$ is a singular point with local group $G_q$, and $f: \mathcal{O}' \rightarrow \mathcal{O}$ is the map with image $\{q\}$.  Then $f$ (or $\mathcal{O}'$ identified with $\{q\}$) is an embedded full suborbifold of $\mathcal{O}$.  By the above definition the normal bundle of $\mathcal{O}'$ in $\mathcal{O}$ is isomorphic to $\mathbb{R}^n //G_q$, which is an orbifold vector bundle over a 0-dimensional (connected) noneffective orbifold with local group $G_q$.

\vspace*{0.2cm}

\noindent {\bf Example} (see Example 10 in \cite{BB15}). Let $\mathcal{Q}$ be the quotient orbifold $\mathbb{R}//\mathbb{Z}_2$, where the  $\mathbb{Z}_2$-action on $\mathbb{R}$ is generated by the reflection with respect to the origin. Let $\Delta: \mathcal{Q} \rightarrow \mathcal{Q} \times \mathcal{Q}$ be the diagonal embedding. Note that $\Delta$
is not a full suborbifold of $\mathcal{Q} \times \mathcal{Q}$. If we still define the normal bundle of the immersion $\Delta$ to be $N_{\Delta}:= \Delta^*T(\mathcal{Q} \times \mathcal{Q})/T\mathcal{Q}$ (but here we view $\mathcal{Q}$ as an effective orbifold), then $N_{\Delta}$ is isomorphic to the tangent bundle $T\mathcal{Q}$, and the only singular point of the total space
of the bundle $N_{\Delta}$ has local group $\mathbb{Z}_2$. But any neighborhood of the image of  $\Delta$ (we also denote this image by $\Delta$) in $\mathcal{Q} \times \mathcal{Q}$ has
a singular point with local group $\mathbb{Z}_2 \times \mathbb{Z}_2$. One sees that the suborbifold $\Delta$   admits no reasonable tubular neighborhood in $\mathcal{Q} \times \mathcal{Q}$.

\vspace*{0.2cm}

Suppose that $\xi$ is a Euclidean orbifold vector bundle over an embedded full suborbifold $\mathcal{O}'$ of $\mathcal{O}$ (so (roughly speaking) it is endowed with a fiberwise inner product structure) with $\text{dim} \hspace*{0.1cm} |E(\xi)|= \text{dim} \hspace*{0.1cm} |\mathcal{O}|$. Let $D(\xi)$ be the (unit) disk bundle of $\xi$. If $\varphi: D(\xi) \rightarrow \mathcal{O}$ is an embedding with the restriction of $\varphi$ to the zero-section of $\xi$ coinciding with the inclusion of $\mathcal{O}'$ in $\mathcal{O}$ (and the corresponding homomorphisms between local groups being isomophisms (but not necessarily the identity maps)), we call $\varphi$ a closed tubular neighborhood of $\mathcal{O}'$  in $\mathcal{O}$. Compare  Section 2 in Chapter VI of \cite{B72},  Section 2 in Chapter III of \cite{K}, and Chapter 7 in \cite{Mu}. Similarly as Theorem 7.1.7 in \cite{Mu}, any embedded full suborbifold $\mathcal{O}'$ in $\mathcal{O}$ (both without boundary) has a closed tubular neighborhood in $\mathcal{O}$. 

Two open tubular neighborhoods $\varphi: E(\xi) \rightarrow \mathcal{O}$ and $\psi: E(\eta) \rightarrow \mathcal{O}$ of an embedded full suborbifold $\mathcal{O}'$ in $\mathcal{O}$ are isotopic if there exist a family of open tubular neighborhoods $\varphi_t: E(\xi) \rightarrow \mathcal{O}$ ($t\in [0,1]$)
of $\mathcal{O}'$ in $\mathcal{O}$ and a  vector bundle isomorphism $\Lambda: E(\xi) \rightarrow E(\eta)$ (over the identity map of $\mathcal{O}'$) such that $\varphi_0=\varphi$, $\varphi_1=\psi \circ \Lambda$,  and the map $[0,1] \times E(\xi)\rightarrow \mathcal{O}$ taking $(t,v) \mapsto \varphi_t(v)$ is smooth. Two closed tubular neighborhoods $\varphi: D(\xi) \rightarrow \mathcal{O}$ and $\psi: D(\eta) \rightarrow \mathcal{O}$ of an embedded full suborbifold $\mathcal{O}'$ in $\mathcal{O}$ are technical isotopic if there exist a family of closed tubular neighborhoods $\varphi_t: D(\xi) \rightarrow \mathcal{O}$ ($t\in \mathbb{R}$)
of $\mathcal{O}'$ in $\mathcal{O}$ and an  orthogonal  bundle map (i.e. a fiberwise isometric map) $\Lambda: D(\xi) \rightarrow D(\eta)$ (over the identity map of $\mathcal{O}'$) such that for some small $\epsilon >0$, $\varphi_t=\varphi$ when $t \leq \epsilon$, $\varphi_t=\psi$ when  $t \geq 1-\epsilon$,  $\varphi_1=\psi \circ \Lambda$,  and the map $\mathbb{R} \times D(\xi)\rightarrow \mathcal{O}$ taking $(t,v) \mapsto \varphi_t(v)$ is smooth. Compare Section 2 in Chapter VI of \cite{B72}  and Chapter 7 in \cite{Mu}.

We need a very simple case of Thom's isotopy extension theorem for orbifolds. Clearly it can be extended to a more general situation.
\begin{lem} \label{lem 2.0}   Let $\mathcal{P}$ be a compact $n$-dimensional orbifold with boundary a manifold,  $\mathcal{O}$ be an $n$-dimensional orbifold (without boundary), and  $f: \mathcal{P}\rightarrow \mathcal{O}$ be an embedding. Suppose that $h: \mathcal{P} \times \mathbb{R} \rightarrow \mathcal{O}$ is a technical isotopy of $f$. Then there exists a diffeotopy  $H: \mathcal{O} \times \mathbb{R} \rightarrow \mathcal{O}$ of $\mathcal{O}$ with compact support such that $h_t=H_t \circ f$ for all $t\in \mathbb{R}$.
\end{lem}
\noindent {\bf Proof}.\ \  The proof can be easily adapted from  that of Theorem 9.5 in  \cite{BJ} and Theorem 7.3.3 in \cite{Mu}  (cf. also p. 179-181 of \cite{M07} and the proof of Theorem 8.6 in \cite{Ka} and Theorem 2.4.2 in \cite{W}).  The desired ambient isotopy comes from the projection of  a suitable global flow on $\mathcal{O}\times \mathbb{R} $, the latter being generated by  a vector field $X$ on $\mathcal{O}\times \mathbb{R}$. Now we proceed to construct $X$. By assumption there exists some small $\epsilon >0$, such that $h_t=f$ when $t \leq \epsilon$, and $h_t=g$ when  $t \geq 1-\epsilon$, where $g: \mathcal{P}\rightarrow \mathcal{O}$ is an embedding. Let $\hat{h}:  \mathcal{P} \times \mathbb{R} \rightarrow \mathcal{O}
 \times \mathbb{R}$ be given by $\hat{h}(x,t)=(h_t(x),t)$. Denote by  $\partial_t$ (resp. $\tilde{\partial}_t$) the vector field on $\mathcal{O} \times \mathbb{R}$ (resp. $\mathcal{P} \times \mathbb{R}$) which projects to 0 on $\mathcal{O}$ (resp. $\mathcal{P}$), and projects to $\frac{\partial}{\partial t}$ on $\mathbb{R}$. On  $ U_1:=\mathcal{O} \times \mathbb{R} \setminus \hat{h}(\mathcal{P} \times [\epsilon, 1-\epsilon])$ let $X_1=\partial_t$. On $U_2:=\hat{h}((\mathcal{P} \setminus \partial \mathcal{P})\times \mathbb{R})$ let $X_2=\hat{h}_*(\tilde{\partial}_t)$. For any point $p$ in $\hat{h}(\partial \mathcal{P}\times \mathbb{R})$, we choose a small open neighborhood $U_p$ of $p$ in $\mathcal{O} \times \mathbb{R}$ and extend the restriction of $\hat{h}_*(\tilde{\partial}_t)$ to the set $U_p \cap \hat{h}(\mathcal{P}\times \mathbb{R})$ to a smooth vector field on $U_p$, denoted by $X_p$, such that the projection of $X_p$ to $\mathbb{R}$ is still $\frac{\partial}{\partial t}$. As usual we choose a partition of unity subordinate to the covering of $\mathcal{O} \times \mathbb{R}$ consisting of $U_1$, $U_2$, and $U_p$'s for $p$ in $\hat{h}(\partial \mathcal{P}\times \mathbb{R})$, and use it to glue $X_1$,  $X_2$, and $X_p$'s, thus get $X$. Clearly $X$ has the desired property. Compare the proof of Theorem 9.5 in  \cite{BJ}, Theorem 7.3.3 in \cite{Mu}, Theorem 2.4.2 in \cite{W}, and Theorem 3.1 in Chapter VI in \cite{B72}.  \hfill{$\Box$}

\vspace*{0.2cm}

As previously mentioned, for any embedded full suborbifold in an orbifold the existence of tubular neighborhoods is known in the existing literature. But I can only show the (ambient) isotopy uniqueness of closed tubular neighborhoods for a smaller class of suborbifolds (i.e. compact perfect suborbifolds (which will be defined below) satisfying an extra technical condition (the latter just for simplicity)).

\vspace*{0.2cm}

\noindent {\bf{Definition}}. Let $\mathcal{O}'$ be an  embedded suborbifold (without boundary, not necessarily connected) of an orbifold $\mathcal{O}$. 
If for any $p\in |\mathcal{O}'|$, there is a good local chart $(U_p, \tilde{U}_p, \Gamma_p, \pi_p)$ of $\mathcal{O}$ around $p$ (see pp. 5790-5791 in \cite{GKRW}) such that $\pi_p^{-1}(U_p \cap |\mathcal{O}'|)$ is a submanifold (of $\tilde{U}_p$) which is pointwise fixed by $\Gamma_p$, we call $\mathcal{O}'$ a perfect suborbifold of $\mathcal{O}$.

\vspace*{0.2cm}
Note that any perfect suborbifold $\mathcal{O}'$ of 
$\mathcal{O}$ is an embedded full suborbifold. Clearly any isolated singular point of an orbifold 
is a perfect suborbifold.

The following result on ambient isotopy uniqueness of  closed tubular neighborhoods of a compact perfect suborbifold (satisfying an extra technical condition) in an  orbifold  will be used extensively in this paper; cf.  p.443 of \cite{BSi} (compare also the Remark on p.312 of \cite{B72} and Theorem 2.9 in \cite{PW}).

\begin{prop} \label{prop 2.1}  Let $\mathcal{O}'$ be a compact perfect suborbifold (without boundary, not necessarily connected) of an orbifold $\mathcal{O}$ with $dim \hspace*{0.1cm}|\mathcal{O}'| <  dim \hspace*{0.1cm} |\mathcal{O}|$. Suppose that  $\varphi: D(\xi) \rightarrow \mathcal{O}$ and $\psi: D(\eta) \rightarrow \mathcal{O}$ are two  closed tubular neighborhoods of $\mathcal{O}'$ in $\mathcal{O}$. Assume  that the boundary of $D(\xi)$ (resp. $D(\eta)$) is a manifold. Then the two closed  tubular neighborhoods are technical isotopic. Moreover, they are ambient isotopic, that is, there exists a diffeotopy $H$  of $\mathcal{O}$ fixing $\mathcal{O}'$ pointwise, with compact support, and with $H_1\circ \varphi=\psi \circ \Lambda$, where $\Lambda: D(\xi) \rightarrow D(\eta)$  is an orthogonal bundle map (over the identity map of $\mathcal{O}'$).
\end{prop}
\noindent {\bf Proof}.\ \   We'll first show that for a perfect suborbifold $\mathcal{O}'$ (not necessarily compact) in $\mathcal{O}$, any two open tubular neighborhoods   $\varphi: E(\xi) \rightarrow \mathcal{O}$ and $\psi: E(\eta) \rightarrow \mathcal{O}$  are isotopic. This can be shown by using local models for the orbifold vector bundles $\xi$  and $\eta$ and the orbifold $\mathcal{O}$, and adapting the proof of Theorem 2.6 in Chapter VI in \cite{B72} (cf. also Theorem 4.6 in \cite{Ka} and Theorem 7.4.6 in \cite{Mu}). We observe that the proof of  Theorem 2.6 in Chapter VI in \cite{B72} can be extended to the weakly equivariant situation when the submanifold $A$ there is contained in the fixed-point set of the $G$-action on $M$; more precisely, the embeddings $\varphi$ and/or $\psi$ there can be only weakly equivariant (that is, equivariant with respect to a group automorphism not necessarily the identity map), moreover, the relevant group automorphisms  need not be the same, if we assume in addition that $A$ is fixed pointwise by $G$. Below we show the extension following the notation  in  the proof of  Theorem 2.6 in Chapter VI in \cite{B72}, but allowing weakly equivariant embeddings in the definition of invariant tubular neighborhoods (compare p. 306 in \cite{B72}); correspondingly, we  extend the notion of isotopy of invariant tubular neighborhoods on  p. 309 of \cite{B72} to the weakly equivariant situation. 

Let $G$ be a compact Lie group, $M$ be a (smooth) $G$-manifold, $A$ be a submanifold of $M$ which is fixed pointwise by $G$. Suppose that $\xi$ and $\eta$ are  (smooth) $G$-vector bundles over $A$ (for the definition see p. 303 in \cite{B72}), and there are weakly $G$-equivariant embeddings $\varphi: E(\xi) \rightarrow M$ and  $\psi: E(\eta) \rightarrow M$, each of whose images is an open neighborhood of $A$ in $M$, such that the restriction of $\varphi$ (resp. $\psi$) to the zero-section of $\xi$ (resp. $\eta$) is the inclusion of $A$ in $M$. Thus  $\varphi(gv)=\alpha(g)\varphi(v)$, and $\psi(gw)=\beta(g)\psi(w)$, for any $g \in G$, $v \in E(\xi)$, and $w\in E(\eta)$, where $\alpha$ and $\beta$ are two automorphisms of the group $G$. W. l. o. g. we may assume that the image of $\varphi$ is contained in the image of $\psi$. Then one can check that the family of embeddings $\varphi_t$ defined by $\varphi_t(v)=\psi(\frac{1}{t}\psi^{-1}\varphi(tv))$ (for $t \neq 0$) satisfies $\varphi_t(gv)=\alpha(g)\varphi_t(v)$ for any  $g \in G$ and $v \in E(\xi)$. That is, each $\varphi_t$ (for $t \neq 0$) is $\alpha$-equivariant. As shown in the proof of  Theorem 2.6 in Chapter VI in \cite{B72} by using Taylor's theorem, for any $v\in E(\xi)$, $\lim_{t \rightarrow 0} \varphi_t(v)$ exists, so one can make sense of $\varphi_0: E(\xi) \rightarrow M$ by defining $\varphi_0(v)=\lim_{t \rightarrow 0} \varphi_t(v)$. Then for any $g \in G$ and $v \in E(\xi)$,
\begin{equation*}
\varphi_0(gv)=\lim_{t \rightarrow 0} \varphi_t(gv)
=\lim_{t \rightarrow 0} \alpha(g)\varphi_t(v)=\alpha(g)\lim_{t \rightarrow 0}\varphi_t(v)=\alpha(g)\varphi_0(v).
\end{equation*}
That is, $\varphi_0$ is also $\alpha$-equivariant. As shown in  \cite{B72}, $\theta:=\psi^{-1}\varphi_0:E(\xi) \rightarrow E(\eta)$ is a vector bundle isomorphism (over the identity map of $A$). One can check that $\theta$ is $\beta^{-1}\alpha$-equivariant. (This also forces us to assume that the submanifold $A$ is fixed pointwise by the compact Lie group $G$.) It follows that the two open invariant (weakly equivariant) tubular neighborhoods $\varphi: E(\xi) \rightarrow M$ and  $\psi: E(\eta) \rightarrow M$ of $A$ in $M$ are weakly equivariantly isotopic. 

Now suppose that both $\xi$ and $\eta$ above are Euclidean $G$-bundles over $A$ (for the definition see p. 304 in \cite{B72}). Let $D(\xi)$ (resp. $D(\eta)$)  be the (unit) disk bundle of $\xi$ (resp. $\eta$). Denote the restriction of  the open weakly equivariant tubular neighborhood $\varphi: E(\xi) \rightarrow M$ (resp. $\psi: E(\eta) \rightarrow M$) to the disk bundle by the same notation $\varphi: D(\xi) \rightarrow M$ (resp. $\psi: D(\eta) \rightarrow M$).
  Consider the polar decomposition of the weakly equivariant vector bundle isomorphism $\theta:  E(\xi) \rightarrow E(\eta)$ (over the identity map of $A$) above,
\begin{equation*}
\theta: E(\xi) \xrightarrow {\theta^{\prime}} E(\xi) 
\xrightarrow {\theta^{\prime \prime}} E(\eta),
\end{equation*}
where $\theta^{\prime}:=\sqrt{\theta^*\theta}$ (with $\theta^*$ being the adjoint of $\theta$ with respect to the $G$-invariant fiberwise inner products on $\xi$ and $\eta$) is self-adjoint and positive definite, and $\theta^{\prime \prime}:= \theta \circ (\theta^{\prime})^{-1}$ is orthogonal  (i.e. a fiberwise isometry). Let $\gamma=\beta^{-1}\alpha$. For any $g \in G$, $x\in A$,  $v\in E(\xi)_x$ and $w\in E(\eta)_x$,   one can compute that 
\begin{equation*}
	\begin{split}
& \langle \theta^*\gamma(g)w, v \rangle =\langle \gamma(g)w, \theta v \rangle =
\langle w,\gamma(g)^{-1}\theta v \rangle \\
 = & \langle w,\gamma(g^{-1})\theta v \rangle = \langle w, \theta g^{-1} v \rangle =  \langle \theta^*w, g^{-1} v \rangle =\langle g\theta^*w, v \rangle.
\end{split}
\end{equation*}
Here we use the assumption that the submanifold $A$  is contained in the fixed-point set of the $G$-action on $M$. It follows that $\theta^*\gamma(g)= g\theta^*$,  $\theta^*\theta g= g\theta^*\theta$, and $\theta^{\prime} g= g\theta^{\prime}$. Then one sees that $\theta^{\prime \prime}$ is $\gamma$-equivariant. Thus as in \cite{B72}, by setting $\theta_t=\theta^{\prime \prime} \circ (t\theta^{\prime}+(1-t)I)$ we see that the two closed invariant (weakly-equivariant) tubular neighborhoods  $\varphi: D(\xi) \rightarrow M$ and $\psi: D(\eta) \rightarrow M$ of $A$ in $M$ are weakly equivariantly technical isotopic; consult pp. 89-90 of \cite{BJ}. 

Then we turn to the orbifold case. Suppose $\varphi: E(\xi) \rightarrow \mathcal{O}$ and $\psi: E(\eta) \rightarrow \mathcal{O}$ are two open tubular neighborhoods of a perfect suborbifold $\mathcal{O}'$ in $\mathcal{O}$. As in the proof of  Theorem 2.6 in Chapter VI in \cite{B72} w. l. o. g. we may assume that the image of $|\varphi|$ is contained in the image of $|\psi|$. In fact we can assume that  $|\varphi|(|E(\xi)|)$ is a sufficiently small open neighborhood of $|\mathcal{O}'|$ in $|\mathcal{O}|$. Note that each of $\varphi$ and $\psi$ may be viewed as a complete orbifold map in the sense of Definition 3.1 in \cite{BB12} and Definition 2.1 in \cite{BB13}.  Using constructions in local models as above (but replace $t$ on the RHS of the equality defining $\varphi_t$ above by $1-t$, that is, set $\tilde{\varphi}_t(v)=\tilde{\psi}(\frac{1}{1-t}
\tilde{\psi}^{-1}\tilde{\varphi}((1-t)v))$ (for $t \neq 1$), where $\tilde{\varphi}$ and $\tilde{\psi}$ are local smooth lifts of $|\varphi|$ and $|\psi|$ respectively, since our definition of isotopy of tubular neighborhoods requires $\varphi_0=\varphi$ instead of $\varphi_1=\varphi$ as in \cite{B72}) we can get from $\varphi$ and $\psi$ a smooth family of complete orbifold maps, denoted by $\varphi_t: E(\xi) \rightarrow \mathcal{O}$ ($t \neq 1$), each of which is an open tubular neighborhood of $\mathcal{O}'$ in $\mathcal{O}$. 

Alternatively, let's check that the $\varphi_t$  ($t \neq 1$) above is well-defined in a certain sense by  using the notion of good map in the sense of Chen-Ruan (see the Appendix in \cite{CR} and Definitions 1.36 and 1.37 in \cite{GMS}); compare Remark 1 in \cite{PR}. Suppose $\text{dim}  \hspace*{0.1cm} |\mathcal{O}|=n$ and  $\text{dim}  \hspace*{0.1cm} |\mathcal{O}^{\prime}|=m$. Let $\varphi: E(\xi) \rightarrow \mathcal{O}$ and $\psi: E(\eta) \rightarrow \mathcal{O}$ be as above. For any $p \in |E(\xi)|$, let  $(U_p, \tilde{U}_p, \Gamma_p, \pi_p)$ be a (good) local chart (local trivialization; cf. Definition 4.1.7 in \cite{CR} and Definition 1.45 in \cite{GMS}) of $E(\xi)$ around $p$ with $\tilde{U}_p=V_p \times \mathbb{R}^{n-m}$ for some open subset $V_p$ of $\mathbb{R}^m$. Given $p^{\prime} \in U_p$, choose a local chart $(U_{p^{\prime}}, \tilde{U}_{p^{\prime}}, \Gamma_{p^{\prime}}, \pi_{p^{\prime}})$  of $E(\xi)$ around $p^{\prime}$ with  $U_{p^{\prime}} \subset U_p \subset |E(\xi)|$ and $\tilde{U}_{p^{\prime}}=V_{p^{\prime}} \times \mathbb{R}^{n-m}$ for some open, simply connected subset $V_{p^{\prime}}$ of $\mathbb{R}^m$. Let $q=|\varphi|(p)$, $q^{\prime}=|\varphi|(p^{\prime})$, $U_q=|\varphi|(U_p)$,  $U_{q^{\prime}}=|\varphi|(U_{p^{\prime}})$, $z=|\psi|^{-1}(q)$, and $z^{\prime}=|\psi|^{-1}(q^{\prime})$. Let
$(U_z, \tilde{U}_z,  \Gamma_z, \pi_z)$ with 
$\tilde{U}_z=V_z \times \mathbb{R}^{n-m}$ for some open subset $V_z$  of $\mathbb{R}^m$, and 
$(U_{z^{\prime}}, \tilde{U}_{z^{\prime}},  \Gamma_{z^{\prime}}, \pi_{z^{\prime}})$ with
$\tilde{U}_{z^{\prime}}=V_{z^{\prime}} \times \mathbb{R}^{n-m}$  for some open, simply connected subset  $V_{z^{\prime}}$ of $\mathbb{R}^m$, be two local charts of $E(\eta)$ around $z$ and $z^{\prime}$ respectively with $U_{z^{\prime}} \subset U_z \subset |E(\eta)|$,  $U_q \subset |\psi|(U_z)=:W_q$ and $U_{q^{\prime}} \subset |\psi|(U_{z^{\prime}})=: W_{q^{\prime}}$. Let $(W_q, \tilde{W}_q, \Gamma_q, \pi_q)$ (resp. $(W_{q^{\prime}}, \tilde{W}_{q^{\prime}}, \Gamma_{q^{\prime}}, \pi_{q^{\prime}})$) be a local chart for $\mathcal{O}$ around $q$ (resp. $q^{\prime}$).  Let $\tilde{U}_q:=\pi_q^{-1}(U_q)$ and  $\tilde{U}_{q^{\prime}}=\pi_{q^{\prime}}^{-1}(U_{q^{\prime}})$.

 Let $\lambda_1:  (\tilde{U}_{p^{\prime}}, \Gamma_{p^{\prime}}, \pi_{p^{\prime}}) \hookrightarrow  (\tilde{U}_p, \Gamma_p, \pi_p)$ be an embedding between local charts of $E(\xi)$ (see Definition 1.2 in \cite{GMS}) which is  fiberwise linear, $\lambda_2:  (\tilde{W}_{q^{\prime}}, \Gamma_{q^{\prime}}, \pi_{q^{\prime}}) \hookrightarrow  (\tilde{W}_q, \Gamma_q, \pi_q)$ be an embedding between  local charts of $\mathcal{O}$, and $\lambda_1^{\prime}:  (\tilde{U}_{z^{\prime}}, \Gamma_{z^{\prime}}, \pi_{z^{\prime}}) \hookrightarrow  (\tilde{U}_z, \Gamma_z, \pi_z)$ be  an embedding between local charts of $E(\eta)$  which is fiberwise linear. The existence of these embeddings can be shown by using the continuation method and definitions of orbifolds and orbifold vector bundles (see for example Chapter 1 in \cite{GMS}) and the simply connectedness condition; compare for example p. 3 in \cite{ALR}, p. 5 in \cite{MP} and p. 467 in \cite{S}.  
 
 Let $\tilde{\varphi}_p: \tilde{U}_p \rightarrow \tilde{U}_q$ be a good  $C^\infty$ lift of the homeomorphism  $|\varphi||_{U_p}: U_p \rightarrow U_q$ (see Definition 1.28 in \cite{GMS}), which is an $\alpha_p$-equivariant diffeomorphism with $|\varphi||_{U_p} \circ \pi_p =\pi_q \circ\tilde{\varphi}_p$, where $\alpha_p: \Gamma_p \rightarrow \Gamma_q$ is an isomorphism;  cf. Lemmas 2 and 3 in \cite{PR}. Let  $\tilde{\varphi}_{p^{\prime}}:=\lambda_2^{-1} \circ \tilde{\varphi}_p \circ \lambda_1$. Note that $\pi_q \circ \tilde{\varphi}_p \circ \lambda_1(\tilde{U}_{p^{\prime}}) =|\varphi| \circ \pi_{p^{\prime}}(\tilde{U}_{p^{\prime}})=U_{q^{\prime}}$, and $\pi_{q^{\prime}}=\pi_q \circ \lambda_2$, so $\pi_{q^{\prime}} \circ \tilde{\varphi}_{p^{\prime}}(\tilde{U}_{p^{\prime}})= \pi_q \circ \tilde{\varphi}_p \circ \lambda_1(\tilde{U}_{p^{\prime}})=U_{q^{\prime}}$, and
 $\tilde{\varphi}_{p^{\prime}}(\tilde{U}_{p^{\prime}})=\tilde{U}_{q^{\prime}}$. Let  $\alpha_{p^{\prime}}:=\lambda_2^{-1} \circ \alpha_p \circ \lambda_1$. Then $\alpha_{p^{\prime}}: \Gamma_{p^{\prime}} \rightarrow \Gamma_{q^{\prime}}$ is an isomorphism.
 Thus 
   $\tilde{\varphi}_{p^{\prime}}: \tilde{U}_{p^{\prime}} \rightarrow \tilde{U}_{q^{\prime}}$ is a good $C^\infty$ lift  of the homeomorphism $|\varphi||_{U_{p^{\prime}}}: U_{p^{\prime}} \rightarrow U_{q^{\prime}} $,  which is an $\alpha_{p^{\prime}}$-equivariant diffeomorphism with $|\varphi||_{U_{p^{\prime}}} \circ \pi_{p^{\prime}} =\pi_{q^{\prime}} \circ\tilde{\varphi}_{p^{\prime}}$, such that  $\tilde{\varphi}_p \circ \lambda_1=\lambda_2 \circ \tilde{\varphi}_{p^{\prime}}: \tilde{U}_{p^{\prime}} \rightarrow \tilde{U}_q$ and  $\alpha_p \circ \lambda_1=\lambda_2 \circ \alpha_{p^{\prime}}: \Gamma_{p^{\prime}} \rightarrow \Gamma_q$; cf. p. 68 of \cite{CR} and Remark 1.32 in \cite{GMS}. 

 Let $\tilde{\psi}_z: \tilde{U}_z \rightarrow \tilde{W}_q$   be a good $C^\infty$ lift of $|\psi||_{U_z}: U_z \rightarrow W_q$, which is a  $\beta_z$-equivariant diffeomorphism with $|\psi||_{U_z} \circ \pi_z =\pi_q \circ\tilde{\psi}_z$, where $\beta_z: \Gamma_z \rightarrow \Gamma_q$  is an isomorphism. 
 Let  $\tilde{\psi}_{z^{\prime}}:=\lambda_2^{-1} \circ \tilde{\psi}_z \circ \lambda_1^{\prime}: \tilde{U}_{z^{\prime}} \rightarrow  \tilde{W}_{q^{\prime}}$, and
 $\beta_{z^{\prime}}:=\lambda_2^{-1} \circ \beta_z \circ \lambda_1^{\prime}: \Gamma_{z^{\prime}} \rightarrow \Gamma_{q^{\prime}}$, we get  
   a good $C^\infty$ lift
 $\tilde{\psi}_{z^{\prime}}$ of $|\psi||_{U_{z^{\prime}}}$, which is a  $\beta_{z^{\prime}}$-equivariant diffeomorphism with $|\psi||_{U_{z^{\prime}}} \circ \pi_{z^{\prime}} =\pi_{q^{\prime}} \circ \tilde{\psi}_{z^{\prime}}$, where $\beta_{z^{\prime}}: \Gamma_{z^{\prime}} \rightarrow \Gamma_{q^{\prime}}$ is an isomorphism, such that $\tilde{\psi}_z \circ \lambda_1^{\prime}=\lambda_2 \circ \tilde{\psi}_{z^{\prime}}: \tilde{U}_{z^{\prime}} \rightarrow \tilde{W}_q$ and   $\beta_z \circ \lambda_1^{\prime}=\lambda_2  \circ \beta_{z^{\prime}}: \Gamma_{z^{\prime}} \rightarrow \Gamma_q$.  

In local charts we define $\tilde{\varphi}_{t,p}$ ($t \neq 1$) by 
\begin{equation*}
\tilde{\varphi}_{t,p}(v)\\:= 
\tilde{\psi}_z(\frac{1}{1-t}\tilde{\psi}_z^{-1} \tilde{\varphi}_p((1-t)v)) 
\end{equation*}
for any $v \in \tilde{U}_p$. One can easily check that $\tilde{\varphi}_{t,p}$  ($t \neq 1$) is $\alpha_p$-equivariant. Now we check that $(\tilde{\varphi}_{t,p}, \alpha_p)$'s are compatible with embeddings between local charts. We already have $\alpha_p \circ \lambda_1=\lambda_2 \circ \alpha_{p^{\prime}}: \Gamma_{p^{\prime}} \rightarrow \Gamma_q$ above. For any $v \in \tilde{U}_{p^{\prime}}$, we compute
\begin{equation*}
	\begin{split}
		& \tilde{\varphi}_{t,p}(\lambda_1(v))\\:= &
		\tilde{\psi}_z(\frac{1}{1-t}\tilde{\psi}_z^{-1} \tilde{\varphi}_p((1-t)\lambda_1(v)))\\= & \tilde{\psi}_z(\frac{1}{1-t}\tilde{\psi}_z^{-1} \tilde{\varphi}_p(\lambda_1((1-t)v)))\\
		= & \tilde{\psi}_z(\frac{1}{1-t}\tilde{\psi}_z^{-1}\lambda_2 \tilde{\varphi}_{p^{\prime}}((1-t)v))\\
		= & \tilde{\psi}_z(\frac{1}{1-t}\lambda_1^{\prime}\tilde{\psi}_{z^{\prime}}^{-1} \tilde{\varphi}_{p^{\prime}}((1-t)v))\\
		= & \tilde{\psi}_z(\lambda_1^{\prime}\frac{1}{1-t}\tilde{\psi}_{z^{\prime}}^{-1} \tilde{\varphi}_{p^{\prime}}((1-t)v))\\
		= & \lambda_2\tilde{\psi}_{z^{\prime}}(\frac{1}{1-t}\tilde{\psi}_{z^{\prime}}^{-1} \tilde{\varphi}_{p^{\prime}}((1-t)v))\\
		=: &  \lambda_2\tilde{\varphi}_{t, p^{\prime}}(v).
	\end{split}
\end{equation*}
It follows that $\tilde{\varphi}_{t,p}$'s ($t \neq 1$) descent to a family of continuous map $|\varphi_t|: |E(\xi)| \rightarrow |\mathcal{O}|$, and we get a family of good $C^\infty$ maps $\varphi_t: E(\xi) \rightarrow \mathcal{O}$ in the sense of Definition 1.37 in \cite{GMS}. Of course 
$\varphi_t$ ($t \neq 1$) is a smooth family of  embeddings.

We can also make sense of $\varphi_1$ as above and show $\theta:=\psi^{-1}\varphi_1: E(\xi) \rightarrow E(\eta)$ is an orbifold vector bundle isomorphism, hence get that the two open tubular neighborhoods $\varphi: E(\xi) \rightarrow \mathcal{O}$ and $\psi: E(\eta) \rightarrow \mathcal{O}$ are isotopic. 

 Finally, let's consider the case of closed tubular neighborhoods of certain special compact suborbifold in orbifolds. Let $\mathcal{O}'$ be a compact perfect suborbifold (without boundary, not necessarily connected) of an orbifold $\mathcal{O}$ with $\text{dim} \hspace*{0.1cm}|\mathcal{O}'| <  \text{dim} \hspace*{0.1cm} |\mathcal{O}|$. Suppose that  $\varphi: D(\xi) \rightarrow \mathcal{O}$ and $\psi: D(\eta) \rightarrow \mathcal{O}$ are two  closed tubular neighborhoods of $\mathcal{O}'$ in $\mathcal{O}$. Assume  that the boundary of $D(\xi)$ (resp. $D(\eta)$) is a manifold. By using an orbifold  version of Lemma 7.4.3 in \cite{Mu} (see also Lemma 2.5.1 in \cite{W}), which can be easily proved by adapting the proof there with the help of Lemma \ref{lem 2.0} above,  the closed tubular neighborhood $\varphi: D(\xi) \rightarrow \mathcal{O}$ (resp. $\psi: D(\eta) \rightarrow \mathcal{O}$) can be extended to an open tubular neighborhood (still denoted by the same notation) $\varphi: E(\xi) \rightarrow \mathcal{O}$ (resp. $\psi: E(\eta) \rightarrow \mathcal{O}$).  So we still have an orbifold vector bundle isomorphism $\theta:  E(\xi) \rightarrow E(\eta)$ (over the identity map of $\mathcal{O}'$) as above. By using the polar decomposition of $\theta$ in local models as above we see that there is a technical isotopy $h: D(\xi) \times \mathbb{R} \rightarrow \mathcal{O}$ between the two embeddings $\varphi: D(\xi) \rightarrow \mathcal{O}$ and $\psi \circ \Lambda:  D(\xi) \rightarrow \mathcal{O}$, where $\Lambda: D(\xi) \rightarrow D(\eta)$  is an orthogonal bundle map (over the identity map of $\mathcal{O}'$); consult pp. 89-90 of \cite{BJ}. 
  Then we can apply Lemma \ref{lem 2.0} and finish the proof.
\hfill{$\Box$}

\vspace*{0.2cm}

\noindent {\bf Remark}.  For the (ambient) isotopy uniqueness of closed tubular neighborhoods of a compact perfect suborbifold $\mathcal{O}'$ in $\mathcal{O}$, the extra condition that the boundary of $D(\xi)$ (and $D(\eta)$) is a manifold should be superfluous, and can be removed if we apply a more general isotopy extension theorem for orbifolds, which can be shown along the lines of the proof of Lemma \ref{lem 2.0} with some extra efforts.

\vspace*{0.2cm}

 Let $D_+$ (resp. $D_-$) be the subset of $\mathbb{S}^{n}//\Gamma$ coming from points in $\mathbb{S}^{n}$ with $x_{n+1}\geq 0$ (resp. $x_{n+1}\leq 0$), and $p_+$ (resp. $p_-$) be the point in $D_+$ (resp. $D_-$) coming from the north (resp. south) pole of $\mathbb{S}^{n}$.

\begin{lem} \label{lem 2.2}  Let $\Gamma$ be as above, and $\mathcal{O}$ be an $n$-dimensional orbifold  diffeomorphic to $\mathbb{S}^{n}//\Gamma$. Write $\mathcal{O}= \Omega_1 \cup  \Omega_2$, where   $ \Omega_2$ is diffeomorphic to $D_-$,  and $\Omega_1 \cap \Omega_2=\partial \Omega_1=\partial \Omega_2$.     Assume   that  $h: \partial D_- \rightarrow \partial \Omega_2$ is a diffeomorphism which extends to a diffeomorphism $\tilde{h}: D_- \rightarrow \Omega_2$. Then  $h$  extends to a diffeomorphism $\tilde{h}': D_+ \rightarrow \Omega_1$.
\end{lem}
\noindent {\bf Proof}.\ \   Let $f: \mathbb{S}^{n}//\Gamma   \rightarrow \mathcal{O}$ be a diffeomorphism. The case that $\Gamma$ is trivial follows from the Cerf-Palais disk theorem, and is  well-known. Now we assume that $\Gamma$ is non-trivial. W.l.o.g. we may assume that $f(p_-)\in \Omega_2$, since otherwise we can replace $f$ by $f$ precomposing with an isometric involution of $\mathbb{S}^{n}//\Gamma$ which exchanges $D_+$ and $D_-$. Then we may view each of $f|_{D_-}: D_- \rightarrow f(D_-)$ and $\tilde{h}: D_- \rightarrow \Omega_2$  as a closed tubular neighborhood of the singular point $f(p_-)=\tilde{h}(p_-)$ in $\mathcal{O}$.   By  Proposition \ref{prop 2.1} there exist a $\Lambda \in \text{Isom}(D_-)$ and a self diffeomorphism $F$ of $\mathcal{O}$  such that  $F\circ f|_{D_-} = \tilde{h} \circ \Lambda$.  Then  $F\circ f|_{D_+}: D_+ \rightarrow \Omega_1$   extends $h\circ \Lambda |_{\mathbb{S}^{n-1}/\Gamma}$.  But $(\Lambda |_{\mathbb{S}^{n-1}/\Gamma})^{-1}$ extends to an isometry of $D_+$, so we are done.
\hfill{$\Box$}

\vspace *{0.2cm}

\noindent {\bf Remark}.  Note that the isometry group $\text{Isom}(\mathbb{S}^{n-1}/\Gamma)$  may have more than two components  when $\Gamma$ is nontrivial.

\vspace *{0.2cm}

   As before, let $B$ be a small, open metric ball in $\mathbb{S}^n// \langle\Gamma, \hat{\sigma} \rangle $   centered at   the image of  the north and  south poles of $\mathbb{S}^n$ under the natural projection $\mathbb{S}^n \rightarrow  \mathbb{S}^n// \langle\Gamma, \hat{\sigma} \rangle $, and $\bar{B}$ be its closure.

\begin{lem} \label{lem 2.3}  Let $\Gamma$,  $\sigma$, and $\hat{\sigma}$ be as before, and $\mathcal{O}$ be an $n$-dimensional orbifold  which is diffeomorphic to $\mathbb{S}^n// \langle\Gamma, \hat{\sigma} \rangle $. Write $\mathcal{O}= \Omega_1 \cup  \Omega_2$, where   $ \Omega_2$ is diffeomorphic to $\bar{B}$,  and $\Omega_1 \cap \Omega_2=\partial \Omega_1=\partial \Omega_2$.     Assume   that  $h: \partial \bar{B} \rightarrow \partial \Omega_2$ is a diffeomorphism which extends to a diffeomorphism $\tilde{h}: \bar{B} \rightarrow \Omega_2$. Then there is a $\lambda \in \text{Isom}(\partial \bar{B})$  such that $h\circ \lambda$  extends to a diffeomorphism $\tilde{h}': \mathbb{S}^n// \langle\Gamma, \hat{\sigma} \rangle \setminus B \rightarrow \Omega_1$.
\end{lem}
\noindent {\bf Proof}.\ \  The case that $\Gamma$ is trivial is treated in \cite{Hu22}. Now we assume that $\Gamma$ is nontrivial. Let $p$ be the (unique) singular point in  $\mathbb{S}^n// \langle\Gamma, \hat{\sigma} \rangle $   coming from the north and  south poles of $\mathbb{S}^n$.  Let $f: \mathbb{S}^n// \langle\Gamma, \hat{\sigma} \rangle   \rightarrow \mathcal{O}$ be a diffeomorphism.  Then we may view both $f|_{\bar{B}}: \bar{B} \rightarrow f(\bar{B})$ and $\tilde{h}: \bar{B} \rightarrow \Omega_2$  as a closed tubular neighborhood of the singular point $f(p)=\tilde{h}(p)$.   By Proposition \ref{prop 2.1} there exist a $\Lambda \in \text{Isom}(\bar{B})$ and a self diffeomorphism $F$ of $\mathcal{O}$  such that  $F\circ f|_{\bar{B}} = \tilde{h} \circ \Lambda$.  Then  $F\circ f|_{\mathbb{S}^n// \langle\Gamma, \hat{\sigma} \rangle \setminus B}: \mathbb{S}^n// \langle\Gamma, \hat{\sigma} \rangle \setminus B \rightarrow \Omega_1$   extends $h\circ \Lambda |_{\partial \bar{B}}$, and we are done.
\hfill{$\Box$}

\vspace *{0.2cm}

Now we start to analyze the structure of orbifold ancient  $\kappa$-solutions which satisfy a certain pinching assumption. (Recall from \cite{B19} and \cite{Hu22} that the orbifold ancient  $\kappa$-solutions we consider are weakly PIC2.) Let $r_1:\mathbb{R}^n \rightarrow \mathbb{R}^n$ be the  reflection w.r.t. the hyperplane $x_1=0$, we also use $r_1$ to denote the restriction of $r_1$ to $\mathbb{S}^{n-1}$. The following result is an improvement of Proposition 2.5 in \cite{Hu22}. (Recall that $C_{\text{PIC}}$, $C_{\text{PIC}1}$ and $C_{\text{PIC}2}$ are the sets of algebraic curvature tensors which have nonnegative isotropic curvature (weakly PIC, see item (i) of Definition 1.1 in \cite{B19}),  weakly PIC1 (see item (ii) of Definition 1.1 in \cite{B19}), and  weakly PIC2 (see item (iii) of Definition 1.1 in \cite{B19}) respectively.)

\begin{prop} \label{prop 2.4} (cf. Corollary 6.7 in \cite{B19},  Theorem 3.4 in \cite{CTZ}, and Proposition 2.5 in \cite{Hu22}) \ \  Let  $n\geq 5$, and  $(\mathcal{O},g(t))$, $t\in (-\infty,T]$, be an
orbifold ancient  $\kappa$-solution of dimension $n$ with at most isolated singularities. Suppose that $(\mathcal{O},g(t))$  satisfies $Rm-\theta R \hspace*{1mm} \text{id} \owedge \text{id} \in C_{\text{PIC}}$ for some uniform constant $\theta>0$, and  there is a spacetime point $(x_0,t_0)$ such that the curvature tensor at $(x_0,t_0)$ lies on the boundary of the cone $C_{\text{PIC}2}$. Then for each $t$, $(\mathcal{O},g(t))$
 is isometric to a shrinking Ricci soliton
$\mathbb{S}^{n-1}/\Gamma \times \mathbb{R}$ or $\mathbb{S}^{n-1}/\Gamma
\times_{\mathbb{Z}_2} \mathbb{R}$. In particular, if $\mathcal{O}$ has no singularities and has exactly one end, it must be diffeomorphic to some $\mathbb{S}^n// \langle\Gamma, \hat{\sigma} \rangle \setminus \bar{B}$,  where $\sigma$ is an isometric involution of the spherical space form $\mathbb{S}^{n-1}/\Gamma$ with no fixed points, and $\hat{\sigma}$ and $B$ are defined as before; if $\mathcal{O}$ has nonempty isolated singularities,  it must be diffeomorphic to $\mathbb{S}^n// (x,\pm x') \setminus \bar{B}$ defined before. Furthermore, there exists a positive constant $\varepsilon_4=\varepsilon_4(n)\leq \varepsilon_1$, where $\varepsilon_1$ is the constant in Proposition A.4 in \cite{Hu22}, with the following property.

\noindent When $\mathcal{O}$ is diffeomorphic to $\mathbb{S}^n// (x,\pm x') \setminus \bar{B}$,  any $2\varepsilon_4$-neck in  $(\mathcal{O}, g(t))$  must be diffeomorphic to $\mathbb{S}^{n-1}\times (0,1)$, and the central cross section $\Sigma$ of any $2\varepsilon_4$-neck in  $\mathcal{O}$ bounds a compact domain $\Omega$; moreover, if $f: \mathbb{S}^{n-1}\rightarrow \Sigma$ is a $2\varepsilon_4$-homothety coming from the $2\varepsilon_4$-neck structure, and $f': \partial (\mathbb{S}^n// (x,\pm x') \setminus B) \rightarrow \mathbb{S}^{n-1}$ is a homothety, either $f\circ f'$ or $f\circ r_1\circ f'$ extends to a  diffeomorphism  $F: \mathbb{S}^n// (x,\pm x') \setminus B \rightarrow \Omega$.

\noindent When $\mathcal{O}$ is  diffeomorphic to some manifold $\mathbb{S}^n// \langle\Gamma, \hat{\sigma} \rangle \setminus \bar{B}$ defined as above,  any $2\varepsilon_4$-neck which is contained in the complement of the homothetically embedded $(\mathbb{S}^{n-1}/ \Gamma)/ \langle \sigma \rangle$ in  $(\mathcal{O}, g(t))$ must be diffeomorphic to $\mathbb{S}^{n-1}/\Gamma \times (0,1)$, and  the central   cross section $\Sigma$ of such a $2\varepsilon_4$-neck  bounds a   compact domain $\Omega$ in $\mathcal{O}$; moreover, if $f: \mathbb{S}^{n-1}/ \Gamma \rightarrow \Sigma$ is a $2\varepsilon_4$-homothety coming from the $2\varepsilon_4$-neck structure, there exist some $\mathbb{S}^n// \langle\Gamma, \hat{\sigma'} \rangle \setminus B'$ and  a homothety $f': \partial \bar{B'} \rightarrow \mathbb{S}^{n-1}/ \Gamma$,  such that   $f\circ f'$ extends to a  diffeomorphism  $F: \mathbb{S}^n// \langle\Gamma, \hat{\sigma'} \rangle \setminus B' \rightarrow \Omega$.
\end{prop}

 \noindent {\bf Proof}.\ \  We only need to prove  the last assertion when $\mathcal{O}$ is diffeomorphic to some manifold $\mathbb{S}^n// \langle\Gamma, \hat{\sigma} \rangle \setminus \bar{B}$  with $\Gamma$  nontrivial, as the remaining part is treated in  Proposition 2.5 in \cite{Hu22}.  In this case  $\mathcal{O}$ is a  smooth manifold with (only) one end, which is diffeomorphic to  $\mathbb{S}^{n-1}/\Gamma \times (0,1)$, and $n$ must be even  as  $\Gamma$  is nontrivial (cf. the proof of Proposition 2.5 in \cite{Hu22}).

As in the proof of Proposition 2.5 in \cite{Hu22}, we can write
$\mathcal{O}=\mathbb{S}^{n-1} \times \mathbb{R}/\widetilde{\Gamma}$, where
$\widetilde{\Gamma}$ is a subgroup of the isometry group of the round cylinder $\mathbb{S}^{n-1} \times
\mathbb{R}$ with $\widetilde{\Gamma}=\Gamma \cup \Gamma^1$, where the second
component of $\Gamma$ (resp. of  $\Gamma^1$) acts on $\mathbb{R}$ as the
identity (resp.  a reflection).  Pick $\tilde{\sigma} \in \Gamma^1$.
Then $\tilde{\sigma} ^2 \in \Gamma$, and $\tilde{\sigma}  \Gamma=\Gamma^1$.  So $\tilde{\sigma} $ induces an isometric involution on $\mathbb{S}^{n-1}/ \Gamma$ (without any fixed points by our assumption on $\mathcal{O}$), which will  be denoted by $\sigma$. The action of $\tilde{\Gamma}$ on $\mathbb{S}^{n-1} \times \mathbb{R}$ will leave exactly one cross section, say $\mathbb{S}^{n-1} \times \{0\}$, invariant. (Compare the beginning of the third paragraph on p. 53 in \cite{CTZ}.) Let $\pi: \mathbb{S}^{n-1} \times \mathbb{R} \rightarrow \mathcal{O}=\mathbb{S}^{n-1} \times \mathbb{R}/\tilde{\Gamma}$ be the natural projection.  Then $\pi( \mathbb{S}^{n-1} \times \{0\})$ is a  homothetically embedded $(\mathbb{S}^{n-1}/ \Gamma)/ \langle \sigma \rangle$ in  $(\mathcal{O}, g(t))$.  Fix any $t$ and write $g=g(t)$. From Proposition  A.4 in \cite{Hu22} we have the following

\vspace *{0.2cm}

\noindent {\bf Claim}.   Assume that $0<\varepsilon_4=\varepsilon_4(n) \leq \varepsilon_1$.  Any $2\varepsilon_4$-neck which is contained in the complement of $\pi( \mathbb{S}^{n-1} \times \{0\})$ in  $(\mathcal{O},g)$ must be diffeomorphic to $\mathbb{S}^{n-1}/\Gamma \times (0,1)$, and the central cross section $\Sigma$ of such a $2\varepsilon_4$-neck bounds a compact, connected, (smooth) submanifold in  $\mathcal{O}$ which contains  $\pi(\mathbb{S}^{n-1} \times \{0\})$.

\vspace *{0.2cm}

Let $N$ be a  $2\varepsilon_4$-neck contained in the complement of $\pi( \mathbb{S}^{n-1} \times \{0\})$ in  $(\mathcal{O},g)$, by the Claim, it must be diffeomorphic to $\mathbb{S}^{n-1}/\Gamma \times (0,1)$. Let  $\psi: \mathbb{S}^{n-1}/\Gamma \times (-(2\varepsilon_4)^{-1}, (2\varepsilon_4)^{-1}) \rightarrow N$ be a diffeomorphism giving the $2\varepsilon_4$-neck structure of $N$, and   $\Sigma=\psi(\mathbb{S}^{n-1}/\Gamma \times \{0\})$.  By the Claim  again, $\Sigma$ bounds a
compact domain, say $\Omega$,  in  $\mathcal{O}$  such that $\Omega$ contains  $\pi(\mathbb{S}^{n-1} \times \{0\})$. Let   $f: \mathbb{S}^{n-1}/\Gamma \rightarrow \Sigma$ be a $2\varepsilon_4$-homothety coming from the $2\varepsilon_4$-neck structure, that is,   $f$ is the same as   $\psi|_{\mathbb{S}^{n-1}/\Gamma \times \{0\}}: \mathbb{S}^{n-1}/\Gamma \times \{0\} \rightarrow \Sigma $ after the canonical  identification of $\mathbb{S}^{n-1}/\Gamma \times \{0\}$ with $\mathbb{S}^{n-1}/\Gamma$. As in the proofs of Proposition 6.17 in \cite{B19}  and of Proposition 2.5 in \cite{Hu22},  we will do metric surgery on $(\mathcal{O},g)$ along the $2\varepsilon_4$-neck $N$ to get a compact Riemannian orbifold $(\mathcal{O}',g')$ which is  weakly PIC2 and strictly PIC if $\varepsilon_4=\varepsilon_4(n)$ is sufficiently small; moreover, there is a point in  $\mathcal{O}'$ where the Ricci curvature (even sectional curvature) is positive.
 We can write $\mathcal{O}'=D_1\cup \Omega$, where $D_1$ is diffeomorphic to $D^n//\Gamma$, and
 $D_1\cap \Omega=\partial D_1=\partial \Omega =\Sigma$.    By  inspecting the surgery procedure we see that $f$ extends to a diffeomorphism   from $ D^n//\Gamma$ to $D_1$.

  As in the proof of Proposition 2.5 in \cite{Hu22}, by running the Ricci flow  we see that $ (\mathcal{O}', g')$ is diffeomorphic to a spherical orbifold with  only one singular point, which has local group $\Gamma$, (and the complement of an open tubular neighborhood of the singular point is diffeomorphic to $\mathbb{S}^n// \langle\Gamma, \hat{\sigma} \rangle \setminus B$,) so the spherical orbifold must be of the form $\mathbb{S}^n// \langle\Gamma, \hat{\sigma'} \rangle$ by using  Lemma 5.2 in \cite{CTZ} and
    Corollary 2.4  in Chapter VI of \cite{B72}, where $\sigma'$ is an isometric involution of the spherical space form $\mathbb{S}^{n-1}/\Gamma$ with no fixed points, and the definition of  $\hat{\sigma'}$ is similar to that of $\hat{\sigma}$  above.  Now we have a diffeomorphism
$H: \mathbb{S}^n// \langle\Gamma, \hat{\sigma'} \rangle=\bar{B'} \cup  (\mathbb{S}^n// \langle\Gamma, \hat{\sigma'} \rangle \setminus B') \rightarrow \mathcal{O}'=D_1 \cup \Omega$, where $B'$ is a small, open metric ball  centered at  the (unique) singular point (coming from the north and  south poles of $\mathbb{S}^n$).  Let $f_1: \partial \bar{B'} \rightarrow \mathbb{S}^{n-1}/ \Gamma$ be a homothety. Note that  $f\circ f_1$  extends to a diffeomorphism  from $\bar{B'}$  to $D_1$.  Then the desired result  follows from Lemma \ref{lem 2.3}.
  \hfill{$\Box$

 \vspace*{0.2cm}
\noindent{\bf Remark}. It seems to be possible to choose $\sigma'=\sigma$ in the above proposition with more efforts, maybe via  torsion invariants (cf. \cite{dR} and \cite{R}); but we don't need this in this paper.
 \vspace*{0.2cm}

The following result  extends    Proposition 6.17 in \cite{B19} and Proposition 2.6 in \cite{Hu22}.

\begin{prop}\label{prop 2.5} (cf. Proposition 6.17 in \cite{B19}  and Proposition 2.6 in \cite{Hu22}) Let  $n\geq 5$ and $\varepsilon_4=\varepsilon_4(n)$ be as in Proposition \ref{prop 2.4}. Let  $(\mathcal{O}, g)$ be a complete, noncompact orbifold diffeomorphic to $\mathbb{R}^n//\Gamma$  with strictly PIC and weakly PIC2. Assume that $N$ is a $2\varepsilon_4$-neck  in $\mathcal{O}$  diffeomorphic to $\mathbb{S}^{n-1}/\Gamma \times (0,1)$, and $\Sigma$ is a central cross section of $N$. When $\Gamma$ is nontrivial we assume that   $\Sigma$ bounds a compact domain $\Omega$ in $ \mathcal{O}$ whose interior contains the unique singular point of $\mathcal{O}$.  Then if $f: \mathbb{S}^{n-1}/\Gamma\rightarrow \Sigma$ is a $2\varepsilon_4$-homothety coming from the $2\varepsilon_4$-neck structure of $N$,  $f$   extends to a  diffeomorphism  $F: D^n//\Gamma  \rightarrow \Omega $.
\end{prop}

\noindent {\bf Proof}. \ \ The case that $\Gamma$ is trivial is treated in Proposition 2.6 in \cite{Hu22}.   When $\Gamma$ is nontrivial, we can  use a surgery  argument  as in the proof of Proposition \ref{prop 2.4}, with the help of Lemma \ref{lem 2.2}.
 \hfill{$\Box$

 \vspace*{0.2cm}

\noindent{\bf Remark}. It seems that in Proposition \ref{prop 2.5} the condition ``When $\Gamma$ is nontrivial we assume that   $\Sigma$ bounds a compact domain $\Omega$ in $ \mathcal{O}$ whose interior contains the unique singular point of $\mathcal{O}$" is redundant, that is, when $\Gamma$ is nontrivial,  $\Sigma$ should always bound a compact domain $\Omega$ in $ \mathcal{O}$ whose interior contains the unique singular point of $\mathcal{O}$. Anyway, when we apply Proposition \ref{prop 2.5} in this paper (see the proof of Propositions \ref{prop 2.6} and \ref{prop 2.8} below), this condition can be verified directly by using Proposition A.4 in \cite{Hu22}.

 \vspace*{0.2cm}

 The following definition of $\varepsilon$-caps  extends/strengthens the corresponding definitions in  \cite{B19} and \cite{Hu22}, and is inspired by  \cite{B19}.

\vspace *{0.2cm}

\noindent {\bf Definition}. (cf.  \cite{B19}  and \cite{Hu22})  Let $\varepsilon_0=\varepsilon_0(n)$ be a small positive constant and $0 < \varepsilon < \frac{1}{4}\varepsilon_0$.  Let $(\mathcal{O},g)$  be an $n$-dimensional Riemannian orbifold  with at most isolated singularities.     Given a point $x_0\in \mathcal{O}$, an open subset
  $U$ of $\mathcal{O}$ is an $\varepsilon$-cap centered at $x_0$ if  $U$ is a topological cap and $U\setminus V$ is an $\varepsilon$-neck, where $V$ is a  compact domain in $U$ diffeomorphic to some $D^n//\Gamma$ or  $\mathbb{S}^n// \langle\Gamma, \hat{\sigma} \rangle \setminus B$ or $\mathbb{S}^n// (x,\pm x') \setminus B$ with $\partial V$ a central cross section of an $\varepsilon$-neck contained in $U$ (i.e., $\partial V$ is the image of $\mathbb{S}^{n-1}/\Gamma \times  \{0\}$ under the diffeomorphism associated to the $\varepsilon$-neck) and with $x_0 \in \text{Int}\hspace*{0.5mm} V$, (we will call $V$ a core of the cap $U$,)  and in addition,

  1. when $U$ is diffeomorphic to $\mathbb{R}^n//\Gamma$, if $\Sigma$ is a central  cross section  of  an $\varepsilon_0$-neck in $U\setminus V$ with an $\varepsilon_0$-homothety $f: \mathbb{S}^{n-1}/\Gamma \rightarrow \Sigma$  coming from the $\varepsilon_0$-neck structure,
    $\Sigma$ bounds  a compact domain $\Omega $ in  $U$, and  $f$ extends to a  diffeomorphism  $F: D^n//\Gamma  \rightarrow \Omega $;

 2. when $U$ is diffeomorphic to  $\mathbb{S}^n// \langle\Gamma, \hat{\sigma} \rangle \setminus \bar{B}$,  if  $\Sigma$ is a central  cross section  of  an $\varepsilon_0$-neck in $U\setminus V$ with an $\varepsilon_0$-homothety $f: \mathbb{S}^{n-1}/\Gamma \rightarrow \Sigma$  coming from the $\varepsilon_0$-neck structure,  $\Sigma$ bounds  a compact domain $\Omega $ in  $U$, and  there exist some $\mathbb{S}^n// \langle\Gamma, \hat{\sigma'} \rangle \setminus B'$ and  a homothety $f': \partial \bar{B'} \rightarrow \mathbb{S}^{n-1}/ \Gamma$,  such that   $f\circ f'$ extends to a  diffeomorphism  $F: \mathbb{S}^n// \langle\Gamma, \hat{\sigma'} \rangle \setminus B' \rightarrow \Omega$;

 3. when $U$ is  diffeomorphic to $\mathbb{S}^n// (x,\pm x') \setminus \bar{B}$,  if   $\Sigma$ is a central  cross section  of  an $\varepsilon_0$-neck in $U\setminus V$ with an $\varepsilon_0$-homothety $f: \mathbb{S}^{n-1}\rightarrow \Sigma$  coming from the $\varepsilon_0$-neck structure, $\Sigma$ bounds  a compact domain $\Omega $ in  $U$, and  $f\circ f'$  extends to a  diffeomorphism $F: \mathbb{S}^n// (x,\pm x') \setminus B \rightarrow \Omega$  for some homothety $f': \partial \bar{B} \rightarrow \mathbb{S}^{n-1}$.

\vspace *{0.2cm}

\begin{prop}\label{prop 2.6} (cf.  Theorem 6.18 in \cite{B19}, Theorem 3.9 in \cite{CTZ}, and Proposition 2.7 in \cite{Hu22}) \ \
Given a small positive constant $\varepsilon$  and a constant $\theta >0$, there exist positive constants $C_1=C_1(n, \theta,\varepsilon)$ and $C_2=C_2(n,\theta,\varepsilon)$,
such that given any noncompact orbifold ancient  $\kappa$-solution
$(\mathcal{O},g(t))$ of dimension $n\geq 5$ (with at most isolated singularities) which satisfies $Rm-\theta R \hspace*{1mm} \text{id} \owedge \text{id} \in C_{\text{PIC}}$ and is not locally isometric to an evolving shrinking round cylinder, for each space-time point $(x_0,t_0)$, there
is an open subset $U$ with
$\overline{B(x_0,t_0,C_1^{-1}R(x_0,t_0)^{-\frac{1}{2}})} \subset U \subset B(x_0,t_0, C_1R(x_0,t_0)^{-\frac{1}{2}})$, which falls into
one of the following categories:

(a) $U$ is a   $\varepsilon$-neck diffeomorphic to $\mathbb{S}^{n-1}/\Gamma \times (-1,1)$  centered at $(x_0,t_0)$, or

(b) $U$ is an $\varepsilon$-cap centered at $(x_0,t_0)$ with any cross-section of the $\varepsilon$-neck contained in the end of it diffeomorphic to $\mathbb{S}^{n-1}/\Gamma$

\noindent for some finite subgroup $\Gamma$ of $O(n)$ acting freely on $\mathbb{S}^{n-1}$;
 moreover, the scalar curvature in $U$ at time $t_0$ is (strictly) between $C_2^{-1}R(x_0,t_0)$ and $C_2R(x_0,t_0)$, and
\begin{equation*}
vol_{g(t_0)}(U) > (C_2|\Gamma|)^{-1} R(x_0,t_0)^{-n/2}.
\end{equation*}
\end{prop}
\noindent {\bf Proof}. \ \ The proof is almost the same as that  of Proposition 2.7 in \cite{Hu22}, with the help of Proposition \ref{prop 2.5}.
 \hfill{$\Box$

\vspace*{0.2cm}

Now we have the following description of the canonical neighborhood property of orbifold ancient $\kappa$-solutions  as in Proposition 2.9 in \cite{Hu22}, but with the strengthened definition of $\varepsilon$-caps.

\begin{prop}\label{prop 2.7} (cf. \cite{P2}, Corollaries 6.20 and 6.22 in \cite{B19},  Theorem 3.10 in \cite{CTZ}, Proposition 2.9 in \cite{Hu22}) \ \
Given   a small positive constant $\varepsilon$  and a constant $\theta >0$, there exist positive constants $C_1=C_1(n, \theta,\varepsilon)$, $C_2=C_2(n,\theta,\varepsilon)$ and $C_3=C_3(n,\theta)$ with the following property:
Suppose $(\mathcal{O},g(t))$ is an orbifold ancient $\kappa$-solution of dimension $n\geq 5$ (with at most isolated singularities) which satisfies
$Rm-\theta R \hspace*{1mm} \text{id} \owedge \text{id} \in C_{\text{PIC}}$. Then either $(\mathcal{O},g(t))$ is  compact and strictly PIC2 for any $t$ (hence diffeomorphic to a spherical orbifold), or for each space-time point $(x_0,t_0)$, there
is  an open subset $U$ of $\mathcal{O}$ with
$\overline{B(x_0,t_0,C_1^{-1}R(x_0,t_0)^{-\frac{1}{2}})} \subset U \subset B(x_0,t_0, C_1R(x_0,t_0)^{-\frac{1}{2}})$ and with the values of the scalar curvature in $U$ at time $t_0$ lying (strictly) between $C_2^{-1}R(x_0,t_0)$ and $C_2R(x_0,t_0)$, which falls into
one of the following two categories:

(a) $U$ is a strong  $\varepsilon$-neck diffeomorphic to $\mathbb{S}^{n-1}/\Gamma \times (-1,1)$  centered at $(x_0,t_0)$, or

(b) $U$ is an $\varepsilon$-cap centered at $(x_0,t_0)$ with any cross-section of the $\varepsilon$-neck contained in the end of it diffeomorphic to $\mathbb{S}^{n-1}/\Gamma$

\noindent for some finite subgroup $\Gamma$ of $O(n)$ acting freely on $\mathbb{S}^{n-1}$, and
\begin{equation*}
vol_{g(t_0)}(U) > (C_2|\Gamma|)^{-1} R(x_0,t_0)^{-n/2};
\end{equation*}
moreover, the scalar curvature in $U$ at time $t_0$ satisfies
the derivative estimates
\begin{equation*}
|\nabla R|< C_3 R^{\frac{3}{2}} \hspace*{8mm} \mbox{and} \hspace*{8mm}
|\frac{\partial R}{\partial t}|< C_3 R^2.
\end{equation*}
\end{prop}

\noindent {\bf Proof}. \ \   One can argue as in the proof of Proposition 2.9 in \cite{Hu22}, using Propositions  \ref{prop 2.4} and   \ref{prop 2.6}.
\hfill{$\Box$

\vspace*{0.2cm}

For orbifold standard solutions we have the following result.

\begin{prop}\label{prop 2.8} (cf. \cite{P2}, Corollary A.2 in \cite{CZ}, Corollary 9.3 in \cite{B19}, and Proposition 2.10 in \cite{Hu22}) \ \
Given   a small positive constant $\varepsilon$,  there exist positive constants $C_1'=C_1'(n,\varepsilon)$, $C_2'=C_2'(n,\varepsilon)$ and $C_3'=C_3'(n)$ with the following property: For each space-time point $(x_0,t_0)$ on an orbifold standard solution
$(\mathbb{R}^n//\Gamma,\hat{g}_\Gamma(t))$ (with at most an isolated singularity), there
is  an open subset $U$ of $\mathbb{R}^n//\Gamma$ with
$\overline{B(x_0,t_0,C_1'^{-1}R(x_0,t_0)^{-\frac{1}{2}})} \subset U \subset B(x_0,t_0, C_1'R(x_0,t_0)^{-\frac{1}{2}})$ and with the values of the scalar curvature in $U$ at time $t_0$ lying (strictly) between $C_2'^{-1}R(x_0,t_0)$ and $C_2'R(x_0,t_0)$, which falls into
one of the following two categories:

(a) $U$ is an   $\varepsilon$-neck diffeomorphic to $\mathbb{S}^{n-1}/\Gamma \times (-1,1)$  centered at $(x_0,t_0)$,  and on $U \times [t_0-\min \{R(x_0,t_0)^{-1},t_0\},t_0]$ the solution $\hat{g}_\Gamma(t)$ is, after scaling with $R(x_0,t_0)$ and shifting  $t_0$ to zero, $\varepsilon$-close to the corresponding subset of the evolving round cylinder $\mathbb{S}^{n-1}/\Gamma \times \mathbb{R}$ over the time interval $[-\min \{t_0R(x_0,t_0), 1\}, 0]$ with scalar curvature 1 at the time zero, and $U$ is disjoint from the surgery cap (in particular, $x_0 \notin B(p_\Gamma, 0, \varepsilon^{-1})$, where $p_\Gamma$ denotes the tip of the orbifold standard solution) when $t_0\leq R(x_0,t_0)^{-1}$, or

(b) $U$ is an $\varepsilon$-cap diffeomorphic to $\mathbb{R}^n//\Gamma$ centered at $(x_0,t_0)$,

\noindent  and $vol_{\hat{g}_\Gamma(t_0)}(U) > (C_2'|\Gamma|)^{-1} R(x_0,t_0)^{-n/2}$; moreover, the scalar curvature in $U$ at time $t_0$ satisfies
the derivative estimates
\begin{equation*}
|\nabla R|< C_3' R^{\frac{3}{2}} \hspace*{8mm} \mbox{and} \hspace*{8mm}
|\frac{\partial R}{\partial t}|< C_3' R^2.
\end{equation*}

\end{prop}

\noindent {\bf Proof}. \ \   The proof is almost the same as that of Proposition 2.10 in \cite{Hu22}, with the help of Proposition \ref{prop 2.5}.
\hfill{$\Box$

\section{Existence of $(r, \delta)$-surgical solutions and applications}

Now we can strengthen Proposition 3.2 in \cite{Hu22}, and recognize the topology of a compact orbifold covered by $\varepsilon$-caps and/or $\varepsilon$-necks.

\begin{prop} \label{prop 3.1} (cf. Proposition 3.2 in \cite{Hu22})\ \  Let $0<\varepsilon_0\leq \varepsilon_3$, where the constant $\varepsilon_3$ is as in Proposition A.7 in \cite{Hu22}.  Fix $0 < \varepsilon <  \tilde{\varepsilon}_1(\varepsilon_0)$, where $\tilde{\varepsilon}_1(\cdot)$ is as in Lemma A.6 in \cite{Hu22}.
Let $(\mathcal{O},g)$ be a closed, connected $n$-orbifold with at most isolated singularities. Suppose that  each point of $\mathcal{O}$ is a center of an
$\varepsilon$-neck or an $\varepsilon$-cap.  Then $\mathcal{O}$ is diffeomorphic  to a spherical orbifold or a connected sum of at most two spherical orbifolds with at most isolated singularities. If we assume in addition  that $\mathcal{O}$ is a manifold, then either $\mathcal{O}$ is diffeomorhic to a spherical space form, or $\mathcal{O}$ is  diffeomorphic  to   a  quotient manifold of $\mathbb{S}^{n-1} \times \mathbb{R}$ by standard isometries.
\end{prop}

\noindent {\bf Proof}. Let $\mathcal{O}$ satisfy the assumption of our proposition.
  With our strengthened definition of $\varepsilon$-caps, we see that $\mathcal{O}$ is diffeomorphic  to a spherical orbifold or a connected sum of at most two spherical orbifolds with at most isolated singularities; cf. the proof  of  Proposition 3.2 in \cite{Hu22}.  Below we only indicate the cases where we can improve the arguments in \cite{Hu22}.

   If $\mathcal{O}$ contains (at least) an  $\varepsilon$-cap of type  $\mathbb{R}^n//\Gamma$ ($|\Gamma|\geq 2$) but no $\varepsilon$-caps of the other types, by using the definition of $\varepsilon$-caps, Lemma A.6 and Proposition  A.7 in \cite{Hu22},  and Theorems 1.5 and 1.9 in Chapter 8 of \cite{Hi}, we see that it  is diffeomorphic to the double of $D^n//\Gamma$, hence it is diffeomorphic to $\mathbb{S}^n// \Gamma$.

Similarly, if $\mathcal{O}$ contains  an  $\varepsilon$-cap of type  $\mathbb{R}^n//\Gamma$ ($|\Gamma|\geq 2$) and an $\varepsilon$-cap of  type  $\mathbb{S}^n// \langle\Gamma, \hat{\sigma} \rangle \setminus \bar{B}$,  we see that it  is diffeomorphic to some $(\mathbb{S}^n// \langle\Gamma, \hat{\sigma'} \rangle \setminus B')\cup_h D^n//\Gamma$ (the notation is as in Section 2), where $h$ is a homothety between the boundaries.  Hence it is diffeomorphic to  a  connected sum  $\mathbb{S}^n// \langle\Gamma, \hat{\sigma'}\rangle \sharp \mathbb{S}^n//\Gamma$ (where the  connected sum occurs at two singular points).  In fact, it is diffeomorphic to $\mathbb{S}^n// \langle\Gamma, \hat{\sigma'}\rangle$;  cf. Theorem 2.2 in Chapter 8  of \cite{Hi} and Proposition 7.7 in \cite{B}.

If $\mathcal{O}$ contains  (at least)  an  $\varepsilon$-cap of type   $\mathbb{S}^n// \langle\Gamma, \hat{\sigma} \rangle \setminus \bar{B}$ ($|\Gamma|\geq 2$) but no caps of the other types,  we see that it  is a smooth manifold diffeomorphic to some $(\mathbb{S}^n// \langle\Gamma, \hat{\sigma}_1 \rangle \setminus B_1)\cup_h (\mathbb{S}^n// \langle\Gamma, \hat{\sigma}_2 \rangle \setminus B_2)$, where $h$ is a homothety between the boundaries.   Clearly it is diffeomorphic to  a  connected sum   $\mathbb{S}^n// \langle\Gamma, \hat{\sigma}_1\rangle \sharp \mathbb{S}^n// \langle\Gamma, \hat{\sigma}_2\rangle$ (where the  connected sum occurs at two singular points).

  We claim that any manifold which is diffeomorphic to a connected sum

 \noindent $\mathbb{S}^n// \langle\Gamma, \hat{\sigma}_1\rangle \sharp \mathbb{S}^n// \langle\Gamma, \hat{\sigma}_2\rangle$ (here $\Gamma$ may be trivial; when $\Gamma$ is nontrivial we require that  the  connected sum occurs at two singular points)   is  diffeomorphic  to   a  quotient manifold of $\mathbb{S}^{n-1} \times \mathbb{R}$ by standard isometries. The reason is as follows.  Let $f_0: [-n,n]\rightarrow [0,\infty)$ be a function satisfying the following conditions:

  (1) $f_0$ is even;

  (2) $f_0$ is smooth, concave and positive on $(-n,n)$;

  (3) $f_0(s)=\sqrt{(n-1)(n-2)}$ for $s\in [-1,1]$;
  and

  (4) there is some $s_1\in (1, n)$ such that $f_0(s)=\sqrt{(n-1)^2-(s-1)^2}$ for $s\in [s_1,n]$.

 Following  Chapter 12 of \cite{MT}, we argue the existence of  the function $f_0$. Let $h_0(s)=\frac{s-1}{\sqrt{(n-1)^2-(s-1)^2}}$ for $s\in [0,n)$, and $\lambda:(-\infty, \infty)\rightarrow [0,1]$ be a non-decreasing smooth function which is identically 0 on $(-\infty, n-\frac{1}{2}+\nu]$ and identically 1 on $[n-\nu, \infty)$  for some  positive number $\nu < \frac{1}{10}$.
 Note that
 \begin{equation*}
\int_1^n h_0(s)\lambda(s+n-\frac{3}{2})ds > \int_{3/2}^n h_0(s)ds > \sqrt{(n-1)(n-2)},
\end{equation*}
 and
 \begin{equation*}
\int_1^n h_0(s)\lambda(s)ds < \int_{n-\frac{1}{2}}^n h_0(s)ds < \sqrt{(n-1)(n-2)},
\end{equation*}
 so there exists some $s_0\in (0, n-\frac{3}{2})$  such that

 \begin{equation*}
\int_1^n h_0(s)\lambda(s+s_0)ds = \sqrt{(n-1)(n-2)}.
\end{equation*}

 Setting
 \begin{equation*}
f_0(s)=\int_s^n h_0(r)\lambda(r+s_0)dr   \hspace*{1mm}  \text{for}  \hspace*{1mm}  s\in [0,n]  \hspace*{2mm}  \text{and}  \hspace*{2mm} f_0(s)=f_0(-s)   \hspace*{1mm}  \text{for}  \hspace*{1mm}  s\in [-n,0],
\end{equation*}
 one can verify that $f_0$ is as desired.

 Now   let
  \begin{equation*}
\Sigma (f_0)=\{(x_1, x_2, \cdot\cdot\cdot, x_{n+1})\in \mathbb{R}^{n+1} \hspace*{1mm}  | \hspace*{1mm} -n \leq x_{n+1} \leq n, \hspace*{2mm} x_1^2+x_2^2+\cdot\cdot\cdot +x_n^2=f_0(x_{n+1})^2\}
\end{equation*}
 be the revolutionary hypersurface in  $\mathbb{R}^{n+1}$ obtained via rotating the graph of the function $f_0$; compare Chapter 12 in \cite{MT}. Then $\Sigma(f_0)$ has an obvious $O(n)\times \mathbb{Z}_2$ symmetry, and $\Sigma(f_0)$ is $O(n)\times \mathbb{Z}_2$-equivariantly diffeomorphic to the standard sphere $\mathbb{S}^n$. It follows that the spherical orbifold $\mathbb{S}^n// \langle\Gamma, \hat{\sigma}\rangle$  is diffeomorphic to the union of a certain Riemannian manifold-with-boundary isometric to some $\mathbb{S}^{n-1}/\Gamma \times_{\mathbb{Z}_2}  [-1,1]$ (which is the quotient of $\mathbb{S}^{n-1}/\Gamma \times [-1,1]$ by the ${\mathbb{Z}_2}$-action  generated by $\hat{\sigma}$  with $\hat{\sigma}(x,s)=(\sigma(x),-s)$ for $x\in
\mathbb{S}^{n-1}/\Gamma $ and $s\in [-1,1]$) and a ``standard'' orbifold cap (which is the  quotient of a ``standard'' (rotationally symmetric) cap (which has an end isometric to  $\mathbb{S}^{n-1} \times [0, \mu]$ for some small $\mu>0$) by $\Gamma$) along their boundaries via an isometry between the boundaries.
  By combining this observation with  the argument on p. 105 of \cite{BJ} and Proposition \ref{prop 2.1},  we see that any connected sum
  $\mathbb{S}^n// \langle\Gamma, \hat{\sigma}_1\rangle \sharp \mathbb{S}^n// \langle\Gamma, \hat{\sigma}_2\rangle$ (here $\Gamma$ may be trivial; when $\Gamma$ is nontrivial we require that  the  connected sum occurs at two singular points)    is diffeomorphic to the union of two Riemannian manifolds-with-boundary along their boundaries via an isometry between the boundaries, where each of the Riemannian manifold-with-boundary  is isometric to some $\mathbb{S}^{n-1}/\Gamma \times_{\mathbb{Z}_2}  [-1,1]$. (When $\Gamma$ is trivial, we use the Cerf-Palais disk theorem to replace Proposition \ref{prop 2.1}.) Such a union has $\mathbb{S}^{n-1} \times \mathbb{R}$-geometry.   (By the way, similarly, by using Proposition \ref{prop 2.1} and (the proof of) Theorem 2.2 in Chapter 8  of \cite{Hi} (cf. also Proposition 7.7 in \cite{B} and Theorem 10.2 in \cite{Ka}) one can show that any connected sum $\mathbb{S}^n//\Gamma  \sharp \mathbb{S}^n//\Gamma$ (where the  connected sum occurs at two singular points when $|\Gamma|\geq 2$) is diffeomorphic to $\mathbb{S}^n//\Gamma$, and  any  connected sum  $\mathbb{S}^n// \langle\Gamma, \hat{\sigma}\rangle \sharp \mathbb{S}^n//\Gamma$ (where the  connected sum occurs at two singular points when $|\Gamma|\geq 2$)  is diffeomorphic to  $\mathbb{S}^n// \langle\Gamma, \hat{\sigma}\rangle$.)

 Thus the proof of Proposition \ref{prop 3.1} is finished.
   \hfill{$\Box$}

\vspace *{0.2cm}
Now we adapt a definition from  \cite{BBM}.

\vspace*{0.2cm}

  \noindent {\bf Definition  }(compare \cite{BBM}, \cite{Hu22}). \ \  A piecewise $C^1$-smooth
  evolving compact Riemannian $n$-orbifold $\{(\mathcal{O}(t), g(t))\}_{t \in I }$ with positive isotropic curvature and with at most isolated singularities is a
  surgical solution to the Ricci flow if it has the following
  properties.

  i. The equation $\frac{\partial g}{\partial t}=-2 \hspace*{0.4mm} \text{Ric}$ is satisfied
  at all regular times;

  ii. For each singular time $t_0$ there is a  finite collection
  $\mathcal{S}$ of disjoint embedded $\mathbb{S}^{n-1}/\Gamma$'s in $\mathcal{O}(t_0)$, and a Riemannian orbifold $\mathcal{O}'$ such that

  (a) $\mathcal{O}'$ is obtained from  $\mathcal{O}(t_0)$ by removing a suitable  open  neighborhood of each  element of $\mathcal{S}$ and
  gluing in a Riemannian orbifold  diffeomorphic to $D^n//\Gamma$ along each boundary component thus produced which is diffeomorphic to $\mathbb{S}^{n-1}/\Gamma$;

 (b) $\mathcal{O}_+(t_0)$ is a union of some connected components of $\mathcal{O}'$ and
 $g_+(t_0)=g(t_0)$ on $\mathcal{O}_+(t_0)\cap \mathcal{O}(t_0)$;

(c) each component of $\mathcal{O}'\setminus \mathcal{O}_+(t_0)$ is
diffeomorphic to a spherical orbifold, or a connected sum of at most two spherical orbifolds.

\vspace *{0.2cm}

  We fix  $\varepsilon_0=\min \{\varepsilon_3, \varepsilon_4\}$, where the constant $\varepsilon_3$ is as in Proposition A.7 in \cite{Hu22},  and the constant $\varepsilon_4$ is as in Proposition \ref{prop 2.4}.  Choose $0 < \varepsilon <  \frac{1}{4} \min \{\frac{1}{10^4n}, \tilde{\varepsilon}_1(\varepsilon_0)\}$, where $\tilde{\varepsilon}_1(\cdot)$ is as in  Lemma A.6 in \cite{Hu22}.

   Given a compact Riemannian orbifold  $(\mathcal{O},g_0)$  of dimension $n\geq 12$ with positive isotropic curvature and with at most isolated singularities, let $(\mathcal{O},g(t))$, $t\in [0,T)$, be the maximal solution to the (smooth)  Ricci flow starting from $(\mathcal{O},g_0)$. Set $\hat{T}=\frac{n}{2\inf_{x\in \mathcal{O}}R(x,0)}$. Then $T \leq \hat{T}$. By  Theorem 1.2 in \cite{B19} there is a continuous family of closed, convex, $O(n)$-invariant sets $\{\mathcal{F}_t ~|~ t\in [0,\hat{T}]\}$ in the vector space $\mathcal{C}_B(\mathbb{R}^n)$ of algebraic curvature tensors on $\mathbb{R}^n$ such that the curvature tensor  of $(\mathcal{O},g_0)$ lies in $\mathcal{F}_0$, the family $\{\mathcal{F}_t ~|~ t\in [0,\hat{T}]\}$ is invariant under the Hamilton ODE $\frac{d}{dt}\text{Rm}=Q(\text{Rm})$, and
\begin{equation*}
 \mathcal{F}_t \subset \{\text{Rm} \hspace*{1mm}|\hspace*{1mm} \text{Rm}-\theta R \hspace*{1mm} \text{id}\owedge \text{id} \in C_{\text{PIC}}\} \cap  \{\text{Rm} \hspace*{1mm}|\hspace*{1mm} \text{Rm}+f(R)\hspace*{1mm} \text{id}\owedge \text{id} \in C_{\text{PIC}2}\}
\end{equation*}
 for any  $t \in [0,\hat{T}]$, where  $f:[0,\infty)\rightarrow [0,\infty)$ is an increasing concave  function satisfying $\text{lim}_{s\rightarrow \infty}\frac{f(s)}{s}=0$, and $\theta$ is a small positive number.  As in \cite{B19}, note that here $f$  and $\theta$   depend only on the initial data  $(\mathcal{O},g_0)$. By (a version of) Hamilton's maximum principle the curvature tensor of $(\mathcal{O},g(t))$ lies in $\mathcal{F}_t$ for any  $t\in [0,T)$.
       Let $\beta=\beta(\varepsilon)$ be the constant given in Lemma A.8 in \cite{Hu22}, and choose
 \begin{equation*}
 C = \max \{100 \varepsilon^{-1}, C_1(n,\theta,\varepsilon), C_2(n,\theta,\varepsilon), C_3(n,\theta), C_1'(n,\beta\varepsilon), C_2'(n,\beta\varepsilon), C_3'(n)\},
 \end{equation*}
 where  the constants on the RHS are from  Propositions \ref{prop 2.7}  and \ref{prop 2.8}.

We can also define the canonical neighborhood assumption  $(CN)_r$ with $(4\varepsilon, 4C)$-control and $(r, \delta)$-surgical solutions as in \cite{P2}, \cite{BBB+}, \cite{BBM}, and \cite{Hu22}.

\vspace *{0.2cm}

As in \cite{Hu22}, with the help of Brendle's curvature pinching estimates (see  Theorem 1.2 and Corollary 1.3 in \cite{B19}), and Propositions \ref{prop 2.7}, \ref{prop 2.8}, and \ref{prop 3.1}, we can show the following theorem.

\begin{thm} \label{thm 3.2}(cf. Theorem 11.2 in \cite{B19}, Theorem 4.6 in \cite{CTZ},  and Theorem 4.11 in \cite{Hu22}) \ \  Let $(\mathcal{O},g_0)$ be a compact Riemannian orbifold of dimension $n\geq 12$ with at most isolated singularities and with positive isotropic curvature.  Let  $\varepsilon$  and $C$ be chosen as above. Then  we can find  positive numbers $\hat{r}$ and $\hat{\delta}$   such that there exists an $(\hat{r}, \hat{\delta})$-surgical solution (with $(4\varepsilon, 4C)$-control) to the Ricci flow starting with $(\mathcal{O},g_0)$, which  becomes extinct in finite time.
\end{thm}

\noindent {\bf Proof of Theorem \ref{thm 1.7}}.  Theorem \ref{thm 1.7} follows from  Theorem \ref{thm 3.2}, Proposition \ref{prop 3.1},  and Lemma 4.1 in \cite{Hu22}; cf. the proof of  Theorem 1.3 in \cite{Hu22}.   \hfill{$\Box$}

\vspace *{0.2cm}

\noindent {\bf Proof of Theorem \ref{thm 1.1}}. Compare the proof of the main theorems in \cite{CTZ} and in \cite{Hu22}. Let $(M, g_0)$ be a Riemannian manifold  satisfying the assumptions of
Theorem \ref{thm 1.1}. From  Theorem \ref{thm 1.7} we see that  $M$ is diffeomorphic to  a connected sum of a finite number of orbifolds, called components, where  each component  is diffeomorphic to either a spherical orbifold with at most isolated singularities or a connected sum of  at most two spherical orbifolds with at most isolated singularities; furthermore, if a component  is a manifold, by Proposition \ref{prop 3.1}, it must be diffeomorhic to  either   a spherical space form, or   a  quotient manifold of $\mathbb{S}^{n-1} \times \mathbb{R}$ by standard isometries.   We denote all  components (whether being a manifold or not) by  $\mathcal{O}_1, \mathcal{O}_2, \cdot\cdot\cdot, \mathcal{O}_k$.

 First we   undo  the connected sum in the orbifold components which are diffeomorphic to $\mathbb{S}^n// (x,\pm x')  \hspace*{1mm} \sharp \hspace*{1mm} \mathbb{S}^n// (x,\pm x')$ (where the connected sum occurs at two regular points) or $\mathbb{S}^n// (x,\pm x')  \hspace*{1mm} \sharp \hspace*{1mm}  \mathbb{R}P^n $, and  do orbifold connected sum operations among the orbifolds thus gotten and all other orbifold components to  resolve the orbifold singularities of the components, which are created by the Ricci flow surgeries. (If a component is a smooth manifold, we don't need to worry about it at this step.) By using the proof of Proposition \ref{prop 3.1} we see that what we get in this way are some manifolds, denoted by $Y_j$, where each of them is  diffeomorphic to either a connected sum
  $\mathbb{S}^n// \langle\Gamma, \hat{\sigma}_1\rangle \sharp \mathbb{S}^n// \langle\Gamma, \hat{\sigma}_2\rangle$ or a self-connected sum of $\mathbb{S}^n//\Gamma$ at the two singular points. (Note that $\mathbb{S}^n// (x,\pm x')$ is isometric to a certain $\mathbb{S}^n//\Gamma$, where $\Gamma \cong \mathbb{Z}_2 < O(n)$ acts on $\mathbb{S}^{n-1}$ antipodally and then on $\mathbb{S}^n$ by suspension.) By the  last claim in the proof of Proposition \ref{prop 3.1} in the former case the resulting manifold admits  $\mathbb{S}^{n-1} \times \mathbb{R}$-geometry.    By an argument similar to that in the  last claim in the proof of Proposition \ref{prop 3.1}, we see that any self-connected sum of $\mathbb{S}^n//\Gamma$ (where $|\Gamma|\geq 2$) at the two singular points must also admit $\mathbb{S}^{n-1} \times \mathbb{R}$-geometry.  (By the way,  a similar argument with Proposition
  \ref{prop 2.1} replaced by Theorem 3.2 in Chapter 8 of \cite{Hi} (see also Corollary 3.7 in Chapter III of \cite{K} and Proposition 2.1 in \cite{Hu17}) shows that  any self-connected sum of $\mathbb{S}^n$ also admits  $\mathbb{S}^{n-1} \times \mathbb{R}$-geometry; but this is well-known.)

Then we perform  manifold connected sums to undo the Ricci flow surgeries which do not introduce orbifold singularities. So we do connected sum among $Y_j$, and a finite number of spherical space forms and a finite number of manifolds each diffeomorphic to     a  quotient manifold of $\mathbb{S}^{n-1} \times \mathbb{R}$ by standard isometries,  which  are components untouched in the previous step. We also  do manifold connected sums which undo the decomposition   of the orbifolds diffeomorphic to $\mathbb{S}^n// (x,\pm x')  \hspace*{1mm} \sharp \hspace*{1mm} \mathbb{S}^n// (x,\pm x')$ or $\mathbb{S}^n// (x,\pm x')  \hspace*{1mm} \sharp \hspace*{1mm}  \mathbb{R}P^n $ performed in the previous step. Thus we recover the original manifold $M$ as a desired manifold connected sum.  \hfill{$\Box$}

\vspace *{0.2cm}

\noindent {\bf Proof of Corollary \ref{cor 1.2}}.  Let $M$ satisfy the assumptions of Corollary \ref{cor 1.2}. Then $M$ has a connected sum decomposition as in the conclusion of Theorem  \ref{thm 1.1}. Note that $n$ is even and $M$ is orientable by assumption, the only orientable $n$-dimensional spherical space form is the sphere $\mathbb{S}^n$, and a manifold-with-boundary of the form $\mathbb{S}^{n-1}/\Gamma_i \times_{\mathbb{Z}_2}  [-1,1]$ is non-orientable  (recall that for even $n$ the action of any finite subgroup of $O(n)$ on $\mathbb{S}^{n-1}$ must be orientation-preserving if it is free),  hence a manifold which is diffeomorphic to the union of two manifolds-with-boundary both of the form $\mathbb{S}^{n-1}/\Gamma_i \times_{\mathbb{Z}_2}  [-1,1]$ along the boundaries is also non-orientable, so  we can assume that in the connected sum decomposition of $M$ each summand is diffeomorphic to the total space of an orientable $\mathbb{S}^{n-1}/\Gamma_j$-bundle over $\mathbb{S}^1$  with structure group $\text{Isom}(\mathbb{S}^{n-1}/\Gamma_j)$. By comparing the first Betti number we see that there is only one such  summand, and we are done.
\hfill{$\Box$}

\vspace *{0.2cm}

\noindent {\bf Proof of Corollary \ref{cor 1.3}}.
Let $(M,g_0)$ be a compact connected
manifold of dimension $n\geq 12$ with  positive isotropic curvature and with $|\pi_1(M)| < \infty$.  By the Seifert-van Kampen theorem, we see that there is at most one nontrivial summand in the connected sum decomposition of $M$ given by the conclusion of Theorem \ref{thm 1.1}, which must be diffeomorphic to a spherical space form.  \hfill{$\Box$}

\vspace *{0.2cm}

\noindent {\bf Proof of Corollary \ref{cor 1.4}}.  Let $M$ satisfy the assumption of Corollary \ref{cor 1.4}. By Corollary 1.2 in \cite{Hu22}, which  may also be viewed as a corollary of Theorem \ref{thm 1.1} here,  there is a finite cover of  $M$ which is diffeomorphic to $\mathbb{S}^n$, or $\mathbb{S}^{n-1} \times \mathbb{S}^1$, or a connected sum of  a finite number of copies of $\mathbb{S}^{n-1} \times \mathbb{S}^1$. From \cite{Mo} we have an explicit construction of the universal cover, denoted by $\widetilde{M}$,  of a connected sum of a finite number of copies of $\mathbb{S}^{n-1} \times \mathbb{S}^1$. Observe that the image of a sphere $\mathbb{S}^k$ under a continuous map from $\mathbb{S}^k$ to $\widetilde{M}$ is compact  for any $k \in \mathbb{N}$.  By using induction, the Blakers-Massey theorem, and the long exact sequence for the relative homotopy groups of a pair, we see that $\pi_i(\widetilde{M})=0$ for $2\leq i \leq n-2$. Alternatively, once we know the universal cover $\widetilde{M}$, we can compute the homology groups  $H_i(\widetilde{M})=0$ for $2\leq i \leq n-2$ (by using the Mayer-Vietoris sequence), and conclude $\pi_i(\widetilde{M})=0$ for $2\leq i \leq n-2$ by using the Hurewicz theorem.  Now the corollary follows.
\hfill{$\Box$}

\vspace *{0.2cm}

\noindent {\bf Proof of Corollary \ref{cor 1.5}}.  Let $M$ satisfy the assumption of Corollary \ref{cor 1.5}. From the proof of Corollary 1.2 in \cite{Hu22}, $M$ admits a finite sheeted regular cover, denoted by $\tilde{M}$, which falls into two cases: (1) $\tilde{M}$ is orientable, and is diffeomorphic to  $\mathbb{S}^n$, or $\mathbb{S}^{n-1}\times \mathbb{S}^1$, or a connected sum of finitely many copies of $\mathbb{S}^{n-1}\times \mathbb{S}^1$; (2) $\tilde{M}$ is non-orientable, and is diffeomorphic to $\mathbb{S}^{n-1} \tilde{\times} \mathbb{S}^1$ or a  connected sum of finitely many copies of $\mathbb{S}^{n-1} \tilde{\times} \mathbb{S}^1$  and $\mathbb{S}^{n-1} \times \mathbb{S}^1$. In the first case we have $H_i(\tilde{M})=0$ for $2  \leq i \leq n-2$ by using the Mayer-Vietoris sequence, and
\begin{equation*}
H^i(M; \mathbb{R}) \cong H^i(\tilde{M}; \mathbb{R})^G=0,   \hspace*{8mm}  2\leq i \leq n-2,
\end{equation*}
(cf. for example, Example 2.6 in \cite{GZ},) where $G$ is the deck transformation group of the covering map $\tilde{M} \rightarrow M$.  In the second case, we take the orientation double cover, denoted by $\tilde{\tilde{M}}$,  of $\tilde{M}$. Of course, an orientable double cover is a regular cover. Moreover we know that  $\tilde{\tilde{M}}$ is diffeomorphic to  $\mathbb{S}^{n-1}\times \mathbb{S}^1$, or a connected sum of finitely many copies of $\mathbb{S}^{n-1}\times \mathbb{S}^1$. Now we can argue as in case (1) to conclude $H^i(\tilde{M}; \mathbb{R})=0$  and $H^i(M; \mathbb{R})=0$  for $2\leq i \leq n-2$.
\hfill{$\Box$}

\vspace *{0.2cm}

\noindent{\bf Remark}. In a preprint posted on the arXiv in January 2008, S. Gadgil and H. Seshadri pointed out that if the fundamental group of a compact manifold
$M$ of dimension $n\geq 4$  with positive isotropic curvature is virtually free, then  $b_i(M)=0$  for $2\leq i \leq n-2$.
 \vspace*{0.2cm}

As in \cite{CTZ}, we also have a topological classification of compact Riemannian orbifold of dimension $n\geq 12$ with at most isolated singularities and with positive isotropic curvature.
\begin{thm} \label{thm 3.3}(cf. Corollary 5.3 in \cite{CTZ})   \ \  Let $(\mathcal{O},g_0)$ be a compact connected Riemannian orbifold of dimension $n\geq 12$ with at most isolated singularities and with positive isotropic curvature.  Then $\mathcal{O}$ is diffeomorphic to a connected sum of some spherical orbifolds with at most isolated singularities and some compact manifolds admitting $\mathbb{S}^{n-1}\times \mathbb{R}$-geometry, where the connected sum occurs at regular points.
\end{thm}
\noindent {\bf Proof}.  The result follows from Theorem \ref{thm 1.7} and the proof of Theorem \ref{thm 1.1} and Proposition \ref{prop 3.1}.
 \hfill{$\Box$}

\vspace *{0.4cm}

\noindent {\bf Acknowledgements}.  I would like to thank Prof. Francis Bonahon for answering my question on the paper  \cite{BSi} concerning tubular neighborhoods of suborbifolds. I would also like to thank Prof. Xu-an Zhao for helpful conversations. I was partially supported by NSFC (12271040).

\vspace *{0.4cm}


\vspace *{0.4cm}

Laboratory of Mathematics and Complex Systems (Ministry of Education),

School of Mathematical Sciences, Beijing Normal University,

Beijing 100875,  People's Republic of China

 E-mail address: hhuang@bnu.edu.cn

\end{document}